\documentstyle[bezier,12pt,leqno]{article}
\textheight9in \def\DATE{\today}
\newtheorem{theorem}{Theorem}
\newtheorem{definition}[theorem]{Definition}
\newtheorem{convention}[theorem]{Convention}

\newtheorem{example}[theorem]{Example}
\newtheorem{lemma}[theorem]{Lemma}
\newtheorem{sublemma}[theorem]{Sublemma}
\newtheorem{proposition}[theorem]{Proposition}

\newtheorem{conjecture}[theorem]{Conjecture}
\newtheorem{principle}[theorem]{Principle}

\catcode`\@=11

\def\ps@myheadings{\let\@mkboth\@gobbletwo
\def\@oddhead{\ifnum\count0=1 \hfill\else
\rightmark \hfil \rm\thepage\fi}%
\def\@oddfoot{\ifnum\count0=1 \hfill \rm 1 \hfill \else
\hfill\fi}
\def\@evenhead%
{\rm\leftmark\hfil\rm\thepage}%
\def\@evenfoot{}\def\sectionmark##1{}
\def\subsectionmark##1{}}

\def\@begintheorem#1#2{\it \trivlist \item[\hskip
 \labelsep{\bf #1\ #2.}]}
\def\@opargbegintheorem#1#2#3{\it \trivlist\item[\hskip%
 \labelsep{\bf #1\ #2.\ (#3)}]}
\def\@endtheorem{\endtrivlist}

\def\@listI{\leftmargin\leftmargini \parsep 1pt plus 2.5pt
 minus 1pt\topsep 5pt plus 2pt minus 3pt%
 \itemsep 0pt plus 2.5pt minus 1pt}
\let\@listi\@listI
\@listi

\def\@sect#1#2#3#4#5#6[#7]#8{\ifnum #2>\c@secnumdepth%
 \def \@svsec {}\else \refstepcounter {#1}\edef \@svsec%
 {\csname the#1\endcsname. \hskip .1em }\fi \@tempskipa%
 #5\relax \ifdim \@tempskipa >\z@ \begingroup #6\relax%
 \@hangfrom {\hskip #3\relax \@svsec }{\interlinepenalty%
 \@M #8.\par }\endgroup \csname #1mark\endcsname {#7}%
 \addcontentsline {toc}{#1}{\ifnum #2>\c@secnumdepth%
 \else \protect \numberline {\csname the#1\endcsname. }%
 \fi #7}\else \def \@svsechd {#6\hskip #3\@svsec #8.%
 \csname #1mark\endcsname {#7}\addcontentsline {toc}{#1}%
 {\ifnum #2>\c@secnumdepth \else \protect \numberline%
 {\csname the#1\endcsname. }\fi #7}}\fi \@xsect {#5}}

\def\section{\@startsection {section}{1}{\z@ }%
 {-3.5ex plus -1ex minus -.2ex}{2.3ex plus .2ex}{\bf }}

\def\thebibliography#1{%
 \section *{References.\@mkboth {REFERENCES}{REFERENCES}}%
 \list {[\arabic {enumi}]}{\settowidth \labelwidth {[#1]}%
 \leftmargin \labelwidth \advance \leftmargin \labelsep %
 \usecounter {enumi}} \def \newblock %
 {\hskip .11em plus .33em minus -.07em} \sloppy \clubpenalty 4000%
 \widowpenalty 4000 \sfcode`\.=1000\relax}

\def\@maketitle{%
 \newpage \null \vskip 2em
 \begin{center}
{\Large\bf \@title \par }
 \vskip 1.5em
 {\large \lineskip .5em
 \begin {tabular}[t]{c}\@author
 \end{tabular}\par}
 \end{center}
  \vskip .8em}

\def\abstract{%
\if@twocolumn \section *{Abstract}
 \else \small\quotation\noindent{\bf Abstract.}\fi}

\catcode`\@=13

\font\tenbold=msbm10 scaled \magstep1
\font\sevenbold=msbm7 scaled \magstep1
\font\fivebold=msbm5 scaled \magstep1
\newfam\boldfam 
\textfont\boldfam=\tenbold
\scriptfont\boldfam=\sevenbold
\scriptscriptfont\boldfam=\fivebold

\def\epi{{\to \hskip -3mm \to}}
\def\ot{\otimes} \def\bk{{\bf k}}
\def\rada#1#2{#1,\ldots,#2}
\def\prada#1#2{#1 + \cdots + #2}
\def\Rada#1#2#3{#1_{#2},\dots,#1_{#3}}
\def\pa{\partial}
\def\znamenko#1{{(-1)^{#1}}}
\def\doubless#1#2{{
\def\arraystretch{.5}
\begin{array}{c}
\mbox{\scriptsize $\scriptstyle #1$}
\\
\mbox{\scriptsize $\scriptstyle #2$}
\end{array}\def\arraystretch{1}
}}
\def\qed{\hspace*{\fill}
\mbox{\hphantom{mm}\rule{0.25cm}{0.25cm}}\\}

\def\parr{{\partial_{{\calR}}}}
\def\dirlim{{{\mathop{{\rm lim}}\limits_{\longrightarrow}}\hskip 1mm}}
\def\span{{\it Span\/}}
\def\paiso{{\partial_{\hskip .5mm \rm iso}}}
\def\calX{{\cal X}}
\def\calY{{\cal Y}}
\def\calS{{\cal S}}
\def\ot{{\otimes}}
\def\op{{\oplus}}
\def\oA{{\overline A}}
\def\oB{{\overline B}}
\def\Ker{\mbox{\it Ker\/}}
\def\cases#1#2#3#4{
                  \left\{
                         \begin{array}{ll}
                           #1,\ &\mbox{#2}
                           \\
                           #3,\ &\mbox{#4}
                          \end{array}
                   \right.
}

\def\desusp{\downarrow\!}

\def\desusp{\downarrow\!}

\def\Gammax{{\Gamma(\xi,\pa\xi)}}
\def\Lie{\hbox{{$\cal L$}{\it ie\/}}}
\def\calP{{\cal P}} 
\def\calR{{\cal R}} \def\calI{{\cal I}}

\def\lbigbrace{{\left(\rule{0pt}{12pt}\right.}}
\def\rbigbrace{{\left.\rule{0pt}{12pt}\right)}}

\def\Iso{{\cal I}{{\it so}}}

\def\alphaiso{{\alpha_{{\rm iso}}}}
\def\Riso{{\cal R}_{{\rm iso }}}
 \def\M{{\cal M}} \def\bfk{{\bf k}}
  \def\bff{{\bf f}}

\def\S{{\cal S}}  \def\SD{{(S,\underline D)}} \def\CD{\SD}
\def\uD{{\underline D}} \def\R{{\cal R}} \def\uS{{\underline S}}
\def\ainfty{{A_\infty}}
\def\P{{\cal P}} \def\catainfty{{{\tt A(\infty\hskip -4mm\infty)}}} 
\def\calO{{\cal O}} 
\def\catchains{{\tt Chain}}
\def\bfV{{\bf V}}
\def\bfW{{\bf W}}

\def\btb{{\bullet \to \bullet}}
\def\Pbtb{{\calP_\btb}}
\def\Rbtb{{\calR_\btb}}
\def\Ass{\mbox{${\cal A}${\it ss}}}
\def\Comm{\mbox{{$\cal C$}\hskip -.3mm {\it om}}}
\def\cathpalg{{\tt Alg}_{h{\cal P}}}
\def\otexp#1#2{{#1^{\ot #2}}}
\def\ainft#1#2#3{{(#1,#2,#3_2,#3_3,\ldots)}}

\def\bfC{{\bf C}}
\def\bfU{{\bf U}}
\def\oX{{\overline X}}
\def\PU{{\P_S *_{\Gamma(X_D)}\P_\uD}}
\def\vlra{{\hbox{$-\hskip-1mm-\hskip-2mm\longrightarrow$}}}
\def\vlla{{\hbox{$\longleftarrow\hskip-2mm-\hskip-1mm-$}}}

\def\colorop #1(#2;#3){{#1}
   \left(\rule{0pt}{15pt}\right.
         \hskip -3mm \begin{array}{c}
	              #3\\#2
                     \end{array}
         \hskip -3mm \left. 
   \rule{0pt}{15pt} \right)
}
\def\coll#1{{\{#1(n)\}_{n\geq 1}}}
\def\id{1\!\!1}
\def\lra{{\longrightarrow}}
\def\rada#1#2{{#1,\ldots,#2}}
\def\prada#1#2{{#1+\cdots+#2}}
\def\orada#1#2{{#1\otimes\cdots\otimes#2}}
\def\otexp#1#2{{#1}^{\ot #2}}

\def\Rada#1#2#3{#1_{#2},\dots,#1_{#3}}

\def\End{\hbox{${\cal E}\hskip -.1em {\it nd}$}}

\def\calQ{{\cal Q}}

\newcommand{\Vtriangle}[6]{
\setlength{\unitlength}{.6cm}
\begin{picture}(5,3.6)(0,-.1)
\thinlines
\put(0,2.5){\makebox(0,0){$#1$}}
\put(5,2.5){\makebox(0,0){$#2$}}
\put(2.5,0){\makebox(0,0){$#3$}}

\put(1,2.5){\vector(1,0){3}}
\put(0.5,2){\vector(1,-1){1.5}}
\put(3,0.5){\vector(1,1){1.5}}

\put(1,1){\makebox(0,0)[r]{{\scriptsize $#5$}}}
\put(2.5,3){\makebox(0,0)[b]{{\scriptsize $#4$}}}
\put(4,1){\makebox(0,0)[l]{{\scriptsize $#6$}}}
\end{picture}
}

\newcommand{\VVtriangle}[6]{
\setlength{\unitlength}{.7cm}
\begin{picture}(5,3.6)(0,-.1)
\thinlines
\put(0,2.5){\makebox(0,0){$#1$}}
\put(5,2.5){\makebox(0,0){$#2$}}
\put(2.5,0){\makebox(0,0){$#3$}}

\put(1,2.5){\vector(1,0){3}}
\put(0.5,2){\vector(1,-1){1.5}}
\put(4.5,2){\vector(-1,-1){1.5}}

\put(1,1){\makebox(0,0)[r]{{\scriptsize $#5$}}}
\put(2.5,3){\makebox(0,0)[b]{{\scriptsize $#4$}}}
\put(4,1){\makebox(0,0)[l]{{\scriptsize $#6$}}}
\end{picture}
}

\newcommand{\square}[8]{
\setlength{\unitlength}{.6cm}
\begin{picture}(5,3.6)
\thinlines

\put(0,3){\makebox(0,0){$#1$}}
\put(5,3){\makebox(0,0){$#2$}}
\put(0,0){\makebox(0,0){$#3$}}
\put(5,0){\makebox(0,0){$#4$}}

\put(-.5,1.5){\makebox(0,0)[r]{{\scriptsize $#6$}}}
\put(5.5,1.5){\makebox(0,0)[l]{{\scriptsize $#7$}}}
\put(2.5,0.5){\makebox(0,0)[b]{{\scriptsize $#8$}}}
\put(2.5,3.5){\makebox(0,0)[b]{{\scriptsize $#5$}}}

\put(1,0){\vector(1,0){3}}
\put(1,3){\vector(1,0){3}}
\put(0,2.5){\vector(0,-1){2}}
\put(5,2.5){\vector(0,-1){2}}
\end{picture}
}

\newcommand{\squarewithdoted}[9]{
\setlength{\unitlength}{.8cm}
\begin{picture}(5,3.6)
\thinlines

\put(0,3){\makebox(0,0){$#1$}}
\put(5,3){\makebox(0,0){$#2$}}
\put(0,0){\makebox(0,0){$#3$}}
\put(5,0){\makebox(0,0){$#4$}}

\put(-.5,1.5){\makebox(0,0)[r]{\scriptsize $#6$}}
\put(5.5,1.5){\makebox(0,0)[l]{\scriptsize $#7$}}
\put(2.5,0.5){\makebox(0,0)[b]{\scriptsize $#8$}}
\put(2.5,3.5){\makebox(0,0)[b]{\scriptsize $#5$}}

\put(1,0){\vector(1,0){3}}
\put(1,3){\vector(1,0){3}}
\put(0,2.5){\vector(0,-1){2}}
\put(5,2.5){\vector(0,-1){2}}

\multiput(.5,.25)(.1,.06){40}{\makebox(0,0){$\cdot$}}
\put(4.5,2.67){\vector(2,1){.1}}
\put(2.5,1.8){\makebox(0,0)[b]{$#9$}}
\end{picture}
}

\newcommand{\squarewithdotedmensi}[9]{
\setlength{\unitlength}{.6cm}
\begin{picture}(5,3.6)
\thinlines

\put(0,3){\makebox(0,0){$#1$}}
\put(5,3){\makebox(0,0){$#2$}}
\put(0,0){\makebox(0,0){$#3$}}
\put(5,0){\makebox(0,0){$#4$}}

\put(-.5,1.5){\makebox(0,0)[r]{{\scriptsize $#6$}}}
\put(5.5,1.5){\makebox(0,0)[l]{{\scriptsize $#7$}}}
\put(2.5,0.5){\makebox(0,0)[b]{{\scriptsize $#8$}}}
\put(2.5,3.5){\makebox(0,0)[b]{{\scriptsize $#5$}}}

\put(1,0){\vector(1,0){3}}
\put(1,3){\vector(1,0){3}}
\put(0,2.5){\vector(0,-1){2}}
\put(5,2.5){\vector(0,-1){2}}

\multiput(.5,.25)(.2,.12){19}{\makebox(0,0){$\cdot$}}
\put(4.5,2.67){\vector(2,1){.1}}
\put(2.5,1.8){\makebox(0,0)[b]{\scriptsize $#9$}}
\end{picture}
}

\newcommand{\squarewithdotedmodified}[9]{
\setlength{\unitlength}{.7cm}
\begin{picture}(5,3.8)
\thinlines

\put(0,3){\makebox(0,0){$#1$}}
\put(5,3){\makebox(0,0){$#2$}}
\put(0,0){\makebox(0,0){$#3$}}
\put(5,0){\makebox(0,0){$#4$}}

\put(-.5,1.5){\makebox(0,0)[r]{{\scriptsize $#6$}}}
\put(5.5,1.5){\makebox(0,0)[l]{{\scriptsize $#7$}}}
\put(2.5,0.5){\makebox(0,0)[b]{{\scriptsize $#8$}}}
\put(2.5,3.5){\makebox(0,0)[b]{{\scriptsize $#5$}}}

\put(1,0){\vector(1,0){3}}
\put(1,3){\vector(1,0){1.8}}
\put(0,2.5){\vector(0,-1){2}}
\put(5,2.5){\vector(0,-1){2}}

\multiput(.5,.25)(.1,.06){40}{\makebox(0,0){$\cdot$}}
\put(4.5,2.67){\vector(2,1){.1}}
\put(2.5,1.8){\makebox(0,0)[b]{{\scriptsize $#9$}}}
\end{picture}
}

\def\SDR#1#2#3#4{{
{
\unitlength=.7pt
\thinlines
\begin{picture}(110.00,20.0)(0.00,14.00)
\put(49.50,0.00){\makebox(0.00,0.00)[t]{$#4$}}
\put(49.50,40.50){\makebox(0.00,0.00)[b]{$#3$}}
\put(97.00,20.00){\makebox(0.00,0.00)[lc]{$#2$}}
\put(5.00,21.00){\makebox(0.00,0.00)[rc]{$#1$}}
\put(21.0,8.0){\vector(-2,1){10.00}}
\put(75.00,33.00){\vector(2,-1){10.00}}
\bezier{100}(16.0,10.50)(49.50,0.50)(79.50,10.50)
\bezier{100}(16.0,30.50)(49.50,40.50)(79.50,30.50)
\end{picture}}
}}

\def\doubless#1#2{{
\def\arraystretch{.5}
\begin{array}{c}
\mbox{\scriptsize $\scriptstyle #1$}
\\
\mbox{\scriptsize $\scriptstyle #2$}
\end{array}\def\arraystretch{1}
}}

\title{Homotopy Algebras are Homotopy Algebras}

\author{Martin Markl%
\thanks{The author was supported by the
grant GA AV \v CR \#1019804 and the Deutsche Forschungsgemeinschaft.}}

\begin{document}

\maketitle

\bibliographystyle{plain}
\baselineskip 18pt plus 2pt minus 1 pt

\begin{abstract}
We prove that strongly homotopy algebras (such as $A_\infty$,
$C_\infty$, sh{} Lie, $G_\infty$,...) are
homotopy invariant concepts in the category of chain complexes. An
important consequence is a rigorous proof that `strongly homotopy
structures transfer over chain homotopy equivalences.'

\noindent 
{\bf Classification:} 55U35, 55U15, 12H05, 18G55

\vskip 3mm
\noindent 
{\bf Plan of the paper:} \ref{intro}.  
                      Introduction 
\hfill\break\noindent 
\hphantom{{\bf Plan of the paper:\hskip .5mm}}  \ref{Furby}. 
                      Colored operads and diagrams -- p.~\pageref{Furby} 
\hfill\break\noindent 
\hphantom{{\bf Plan of the paper:\hskip .5mm}} \ref{kam_nas_vystehuji?}.
                      Resolutions, homotopy diagrams and extensions
                      -- p.~\pageref{kam_nas_vystehuji?}
\hfill\break\noindent 
\hphantom{{\bf Plan of the paper:\hskip .5mm}}
                                    \ref{zase_starosti}.
                      Structure of strongly homotopy diagrams
                      -- p.~\pageref{zase_starosti}
\hfill\break\noindent 
\hphantom{{\bf Plan of the paper:\hskip .5mm}}
                                    \ref{jeste_jeden_medvidek}.
                       Proofs of the moves -- p.~\pageref{jeste_jeden_medvidek}
\hfill\break\noindent 
\hphantom{{\bf Plan of the paper:\hskip .5mm}}
                                    \ref{jeden_medvidek}.
                       Appendix (remaining proofs) 
                       -- p.~\pageref{jeden_medvidek}
\end{abstract}

\section{Introduction}
\label{intro}

Historically, the first examples of homotopy algebras were
strongly homotopy associative algebras, also called 
$A_\infty$-algebras or sh associative algebras, discovered, as 
structures intrinsically present in the singular chain complex of a
loop space, by Jim Stasheff~\cite{stasheff:TAMS63} in~1963.

A {\em strongly homotopy associative algebra\/} is a chain complex
$V = (V,\pa)$ together with a bilinear product $\mu_2 : V\ot V \to V$
which is
associative only {\em up to a homotopy\/} $\mu_3 : \otexp V3 \to
V$, that is, for each $a,b,c\in V$,
\[
\mu_2(\mu_2(a,b)c) - \mu_2(a,\mu_2(b,c)) = [\mu_3,\pa](a,b,c),
\] 
where 
\[
[\mu_3;\pa](a,b,c) := \mu_3(\pa a,b,c)+(-1)^{|a|}
\mu_3(a,\pa b,c) + (-1)^{|a|+|b|}
\mu_3(a,b, \pa c) + \pa \mu_3(a,b,c)
\]
is the differential induced on
${\it Hom}(\otexp V3,V)$ by $\pa$. Very crucially, there are some higher
coherence relations involving $\mu_2$, $\mu_3$ and higher homotopies
$\mu_n : \otexp Vn \to V$, $n\geq 4$. Namely, for each $n\geq 2$ and
$\Rada a1n \in V$,
\begin{equation}
\label{An}
\sum_{\doubless{i+j=n+1}{i,j \geq 2}}
\sum_{s=0}^{n-j}
\znamenko{\epsilon} \cdot 
\mu_i(\Rada a1s,\mu_j(\Rada a{s+1}{s+j}),\Rada a{s+j+1}n)
=
[\mu_n,\pa](\Rada a1n),
\end{equation}
where $\epsilon := j+s(j+1) + j(|a_1|+ \cdots + |a_s|)$. 
For a gentle introduction to these structures, we
recommend~\cite{markl:shalg}.

Strongly homotopy (sh{})
morphisms of $A_\infty$-algebras were introduced in~1965 by 
Allan Clark~\cite{clark:PJM65} (though strongly homotopy maps of
topological spaces with additional structures were considered 
by M.~Sugawara~\cite{sugawara:hc} in~1960).  
If $\bfV = \ainft V{\pa}{\mu}$ 
and $\bfW = \ainft W{\pa}{\nu}$ are two $A_\infty$-%
algebras, then a {\em morphism\/} ${\bf f}: \bfV \to \bfW$ 
is a {\em sequence\/} ${\bf f}
=(f_1,f_2,f_3,\ldots)$ of multilinear maps $f_n : \otexp Vn \to W$, $n
\geq 1$,
such that $f_1$ is a chain map which commutes with the products
$\mu_2$ and $\nu_2$ up to a homotopy $f_2$,
\begin{equation}
\label{whistle}
f_1(\mu_2(a,b)) - \nu_2(f_1(a),f_1(b)) = [f_2,\pa](a,b),
\end{equation}
and the following higher coherence relations are satisfied, 
for each $n
\geq 2$ and $\Rada a1n \in V$:
\begin{eqnarray} 
\label{Fn}
&&
\sum_{k=2}^{n} \sum_{r_1+\dots+r_k = n}
(-1)^{\eta}\cdot 
\nu_k(f_{r_1}(\Rada a1{r_1}),\dots,f_{r_k}(\Rada a{n-r_k+1}n))- 
\\
\nonumber 
&&- \hskip -3mm \sum_{\doubless{i+j=n+1}{j\geq 2}}\sum_{s=0}^{n-j}
(-1)^{\nu}\cdot
f_i(\Rada a1{s},\mu_j(\Rada a{s+1}{s+j}),\Rada
a{s+j+1}{n}) = [f_n,\pa](\Rada a1n),
\end{eqnarray}
where 
\[
\eta : = \sum_{1\leq i<j \leq k}r_i(r_j+1) + \sum_{1 < i \leq
k}(r_i +1)(|a_1| + \cdots + |a_{r_1+ \cdots + r_{i-1}}|)
\] 
and 
\[
\nu : = n+s(j+1) + j(|a_1| + \cdots + |a_s|).
\] 

The morphism ${\bf f} =(f_1,f_2,f_3,\ldots)$ is called {\em strict\/},
if $f_i = 0$ for $i \geq 2$. Thus a strict morphism is simply a chain
map $f = f_1 : (V,\pa) \to (W,\pa)$ which commutes 
with the $A_\infty$-structure maps, that is, $f \mu_n = \nu_n f^{\otimes n}$ 
for $n \geq 2$.

If we let $\catchains$ to denote the category of chain complexes over
a ring $R$ and $\catainfty$ the category of $A_\infty$-algebras and their sh{}
morphisms, then the correspondence 
\[
\bfV = \ainft V{\pa}\mu  \longmapsto (V,\pa),\
{\bf f} = (f_1,f_2,\ldots) \longmapsto f_1,
\]
defines the `forgetful' functor $\# : \catainfty \to
\catchains$. Observe that the functor $\#$ is {\em not} 
injective on morphisms.

More examples of homotopy algebras appeared later on. 
Let us name at least
strongly homotopy associative commutative
algebras~\cite{kadeishvili:ttmi85,markl:JPAA92}, also called
$C_\infty$ or balanced $A_\infty$-algebras, and strongly
homotopy Lie algebras, also called sh Lie or
$L_\infty$-algebras~\cite{lada-stasheff:IJTP93,lada-markl:CommAlg95}. 
Then, in the early 90's, with the renaissance of operads, the 
heavens opened wide and examples of homotopy algebras became abundant.

With hindsight, one might say that by a (strongly) homotopy $\P$-algebra, where
$\P$ is a dg-operad (i.e.~an operad in the category of chain complexes), 
is today understood an algebra over a {\em cofibrant model} of $\P$. A~cofibrant 
model of $\P$ is, by definition, a free operad
$\Gamma(E)$ equipped with a differential $\pa$ and a map $\alpha :
(\Gamma(E),\pa)\to \P$ that is a homology isomorphism. We also assume the
existence of a suitable filtration on the space of generators $E$,
see Definition~\ref{cof} for details.

This definition is broad enough to comprise all objects which
would probably be called homotopy algebras, including various
types of homotopy Gerstenhaber algebras arising in connection with
quantum field theory~\cite{getzler-jones:preprint}, the
Deligne conjecture~\cite{voronov:99} 
and, later, in the proof of Kontsevich's formality, 
see~\cite{kontsevich:defquant,kontsevich:motives}.

If $R = \bk$, the field of characteristic zero, 
not only does there exist a cofibrant
model for any operad $\P$, but there even exists, under the assumption
$\calP(0) =0$ and $\P(1)=\bfk$,
a~very special type of this model, 
called the {\em minimal model\/} of $\P$.
By definition, the minimal model of $\P$ is a cofibrant model
$\M_{\P}= (\Gamma(E),\pa)$ of $\P$ as above such that
$\pa(E)$ consists of decomposable elements of the free operad
$\Gamma(E)$. 

In case $\P$ is quadratic Koszul, the minimal model of $\P$
is isomorphic to the
dual bar construction on the quadratic dual $\P^!$ and our conception
coincides with 
that of~\cite{ginzburg-kapranov:DMJ94}. In particular, if $\P
=\Ass$, the operad for associative algebras, strongly homotopy $\P$-algebras
are precisely $A_\infty$-algebras discussed above. 
In more detail, the operad $\Ass$ for associative algebras is given
by
\[
\Ass{} := \Gamma(\mu)/(\mu(\mu \otimes \id) - \mu(\id \otimes \mu)),
\]
the free operad on one bilinear operation $\mu$, modulo the
associativity. The minimal model for the operad $\Ass$ is
\[
\Ass{} \stackrel{\alpha}{\longleftarrow}
\left(
\Gamma(
\mu_2,\mu_3,\mu_4,\ldots),\pa
\right),
\]
where $\alpha(\mu_2) = \mu$
while $\alpha$ is trivial on the remaining generators.
The differential is defined by 
\[
\pa (\mu_n)
=
\sum_{\doubless{i+j=n+1}{i,j \geq 2}} \sum_{s=0}^{n-j}
(-1)^{j+s(j+1)} \mu_{i}(\id^{\ot s} \ot \mu_j \ot
\id^{\ot i-s-1}),
\]  
see~\cite[Example~4.8]{markl:zebrulka}.
Comparing it with~(\ref{An}), we immediately see that $\ainfty$-%
algebras are algebras over the operad $(\Gamma(
\mu_2,\mu_3,\mu_4,\ldots),\pa)$ defined above.

Let us make another brief excursion to topology. With a great
simplification, we might say that key problems in
algebraic topology of the sixties were related to
homotopy invariant structures on topological
spaces. For instance, given an
associative topological monoid $Y$, is it true that each space $Y'$
of the same homotopy type inherits from $Y$ an induced associative monoid
structure? 

The answer is no, but to understand what of the monoidal structure
remains preserved under a change of the homotopy type of the underlying space
required a considerable
effort. In this concrete case, the answer was that the corresponding
homotopy invariant structure is that of a
`strongly homotopy associative monoid,' which is, by definition, a
space acted on by the topological
operad ${\cal K} =\{K_n\}_{n\geq 1}$ of the 
associahedra~\cite{stasheff:TAMS63}.

A conceptual solution of the problem was 
provided by Boardman \& Vogt and
their $W$-construction~\cite{boardman-vogt:73}. 
In today's language, one would say that
the $W$-construction gives, to any topological
operad $\calO$, its specific cofibrant resolution $W\calO$. 
Homotopy invariant algebras are then algebras over the operad
$W\calO$. Another, in a sense dual, solution
of the problem based on `lax' actions of operads was given by
Lada~\cite{lada:1976}.

Let us come back to the realm of 
algebra and try to formulate what homotopy
invariance might mean here. 
The underlying category is the category $\catchains$ of chain complexes
and their chain maps, and homotopy is
provided by chain homotopy of maps.  
If we assume that the ground ring is a field,
homotopy boils down to homology. Thus, for instance, two maps
$f,g: (V,\pa) \to (W,\pa)$ are (chain)
homotopic if and only if $H_*(f) = H_*(g)$.

Suppose we have a `category' (in a suitable sense, see below)
of strongly homotopy $\P$-algebras $\cathpalg$,
with a forgetful functor $\# : \cathpalg \to \catchains$. For an sh{}
morphism ${\bf f}$ of strongly
homotopy $\P$-algebras, we call the chain map
${\bf f}_{\#}$
the {\em underlying map\/} of ${\bf f}$. If $g$ is a chain map and ${\bf f}$
such that ${\bf f}_{\#} = g$, we say that ${\bf f}$ is a {\em strongly
homotopy
$\P$-structure\/} on the map $g$. Let us remark that 
there are certain problems with the categorial structure of
strongly homotopy 
$\P$-algebras (see~\cite{boardman-vogt:73} where, for topological
spaces, these subtleties required the use of weak Kan categories). 
This delicacy, discussed also in~\cite{markl:ho}, 
will however play no r\^ole in this paper.

What homotopy invariance means for topological spaces was formulated 
by Boardman and Vogt 
in the introduction 
to~\cite{boardman-vogt:73}. Their conditions, which we will also call 
{\em moves\/}, reformulated to algebra read as follows.

\begin{enumerate}
\item[{\bf (M1)}]
For each strongly homotopy $\P$-algebra $\bfV$, chain complex
$W = (W,\partial)$ and a chain homotopy equivalence 
$f: (V,\pa) \to (W,\pa)$, there exist a sh{}
$\P$-structure $\bfW$ on the chain complex
$(W,\pa)$ and a strongly homotopy $\P$-structure ${\bf f}:\bfV \to \bfW$
on the map $f$.
\item[{\bf (M2)}]
Suppose $\bfV$ and $\bfW$ are two strongly
homotopy $\P$-algebras and ${\bf f}:\bfV \to
\bfW$ a sh{} $\P$-algebra morphism. Suppose that $g:(V,\pa) \to (W,\pa)$
is a chain map that is chain homotopic to
$\bff_{\#}$. Then there exists a strongly
homotopy $\P$-structure on $g$.
\item[{\bf (M3)}]
Suppose that ${\bf f}:\bfV\to
\bfW$ is a strongly 
homotopy $\P$-algebra morphism such that $\bff_{\#}$ is a chain homotopy
equivalence. 
Suppose that  $g:(W,\pa) \to (V,\pa)$ is a chain homotopy
inverse of $\bff_{\#}$. Then there exists a strongly
homotopy $\P$-structure on $g$.
\end{enumerate}

Moves {\bf (M1)}~and {\bf (M3)} clearly imply: 
\begin{itemize}
\item[{\bf (M1')}]
For each strongly homotopy $\P$-algebra ${\bf W}$, 
chain complex
$V = (V,\partial)$ and a chain homotopy equivalence 
$f: (V,\pa) \to (W,\pa)$, there exist a strongly homotopy 
$\P$-structure ${\bf V}$ on
$(V,\pa)$ and a strongly homotopy $\P$-structure ${\bf f}:{\bf V} \to
{\bf W}$ on the map $f$. 
\end{itemize}

Let us add, for the completeness of the exposition, also the
following move which says that strongly homotopy algebras and
their strongly homotopy maps form a certain weak category.

\begin{itemize}
\item[{\bf (M4)}]
For each sequence ${\bf V}_1 \stackrel{{\bf f}_1}{\lra} {\bf V}_2
\stackrel{{\bf f}_2}{\lra} 
\cdots
\stackrel{{\bf f}_{n-1}}{\lra} 
 {\bf V}_{n} 
\stackrel{{\bf f}_n}{\lra} {\bf V}_{n+1}$ 
of sh{} $\P$-algebra
morphisms and for each chain map $g : (V_1,\pa) \to (V_{n+1},\pa)$,
homotopic to the composition $\#({\bf f}_n) \cdots \#({\bf f}_1)$, there exist
a strongly homotopy $\P$-structure on the map $g$.
\end{itemize}

The aim of this paper is to prove the above `moves' for the case
when the ground ring ${\bf k}$ is a field of 
characteristic zero -- see Section~\ref{jeste_jeden_medvidek}.
We show that
they are, in fact, special cases of a certain
homotopy extension problem for diagrams.
Our main theorem then states that, roughly speaking, 
the homotopy extension problem over cofibrant
operads can be solved if and only if
the corresponding `classical' extension
problem has a solution.

This is our understanding of the following principle that has
been latent in topology since~\cite{boardman-vogt:73}:

\begin{principle}
Algebras over a cofibrant operad are homotopy invariant.
\end{principle}

The most important seems to be move {\bf (M1)} and
its `inverse' {\bf (M1')} which says that `strongly homotopy
structures can be transfered over chain homotopy equivalences.'
Various special cases of this statement were proved by many authors
(Kadeishvili~\cite{kadeishvili:ttmi85},
Lambe, Gugenheim, 
Stasheff~\cite{lambe-stasheff:mm87,gugenheim-stasheff:BSMB86}, 
Huebschmann-Kadeishvili~\cite{huebschmann-kadeishvili:MZ91},
Hess~\cite{hess}, Johannson-Lambe~\cite{johansson-lambe:96},...) 
and used as
a folklore in deformation quantization papers of 
Kontsevich~\cite{kontsevich:defquant,kontsevich:motives}.
Our theory explains the nature of `side conditions' for strong
deformation data in perturbation theory and shows why they are
important for the existence of the `transferred structure,' see
Proposition~\ref{why_side_conditions}. This direction is further
pursued in~\cite{markl:ip}.

\noindent
{\bf What do the moves say?}
{}For a given operad $\calP$, each of
the above moves can be translated to a very concrete statement about
very concrete algebraic structures.

As we already noticed, each Koszul quadratic operad
has a very explicit minimal model. The same is true also for the
minimal model of the colored operad $\calP_{\btb}$ 
describing homomorphisms of $\calP$-algebras discussed
in Section~\ref{Furby}, which gives rise to strongly
homotopy morphisms of strongly homotopy $\calP$-algebras. For example,
${\bf (M1)}$ with $\calP = \Ass$ says:

{\em
For each $A_\infty$-algebra ${\bf V} = (V,\pa,\mu_2,\mu_3,\ldots)$, each
chain complex $W = (W,\pa)$ and a chain map $f :(V,\pa) \to (W,\pa)$
that is a chain homotopy equivalence, there exists an
$A_\infty$-structure ${\bf W} = (W,\pa,\nu_2,\nu_3,\ldots)$ on
$(W,\pa)$ and a strongly homotopy morphism ${\bf f} = (f_1,f_2,\ldots)
: {\bf V} \to {\bf W}$ such that $f_1 = f$.
}
 
We leave as an exercise to interpret the remaining moves. 
A similar explicit translation can be made also for other Koszul
quadratic operads as $\calP = \Lie$, 
the operad for Lie algebras, in which case the
appropriate notions of strongly homotopy Lie algebras and their strongly
homotopy morphisms have been worked out in
~\cite{lada-markl:CommAlg95}. 
For general operads, results of~\cite{markl:ho} can be used.

\noindent
{\bf Comparison with other approaches.}
There have recently appeared other formulations of homotopy invariance
of strongly homotopy algebras 
(see~\cite{berger-moerdijk:02,hinich:preprint00}). 
Instead of analyzing homotopy diagrams, they use various forms 
of `weak equivalences' of strongly homotopy structures. 

For instance, in~\cite{hinich:preprint00} two sh{} algebras are weakly
equivalent, if they are connected by a chain of {\em strict\/}
homomorphisms (not necessarily pointing in the same direction)
that are homology isomorphisms. 
We think that the advantage of our approach is that our
notion of `weak equivalence' can be made, 
for each particular application, very explicit.

\noindent
{\bf Conventions.}
Throughout the paper, we work with {\em algebras without units\/}. In most
situations, the results can be easily translated to the unital case,
but there are important exceptions, such as move {\bf (S)} of
Example~\ref{Salvatore}. 

All algebraic objects are defined over a {\em field ${\bf k}$ 
of characteristic zero\/}. Working over a field simplifies the use of various
forms of K\"unneth formulas and also implies that all
finite-dimensional modules are free. Characteristic zero
assumption is then used in constructions of (minimal) bigraded 
models where we need to `average' over symmetric 
groups~\cite[p.~1476]{markl:zebrulka}. We believe that 
some results of this paper are true also over more general
rings, but from the point of applications the characteristic zero
assumption is very reasonable and saves
us from many potential troubles.

By an operad we always mean an operad in the symmetric monoidal 
category of chain complexes over ${\bf k}$. When we speak about
strongly homotopy $\calP$-algebras, we always assume that $\calP$ has
{\em trivial differential, $\calP(0) = 0$ and $\calP(1) \cong \bfk$\/}.

\noindent
{\bf Acknowledgment.} I would like to express my gratitude to
Rainer~Vogt for his hospitality and numerous discussions 
during my multiple visits of
Osnabr\"uck. My thanks are
also due to Jim Stasheff
for careful reading the manuscript and many
useful remarks, to Paolo Salvatore who pointed the `move' of
Example~\ref{Salvatore} to me and taught me a lot about closed model
categories and to Vladimir Hinich who explained to me the difference
between $B_\infty$ and $G_\infty$. Also the referee's suggestions were
extremely helpful.

\section{Colored operads and diagrams}
\label{Furby}

Since we will be primarily interested in algebras and their {\em
diagrams\/}, we need to introduce (multi)colored operads as a tool 
describing these objects. These colored operads are, in fact,
equivalent to the classical `colored {\small PROP}s' 
of~\cite{boardman-vogt:73}.
In~\cite{baez-dolan:III} they are called `many-sorted' operads.

Fix a (finite) {\em set of colors\/} ${\bf C}= \{\rada {c_1}{c_C}\}$ and
consider an operad $\P = \coll \P$ such
that each $\P(n)$ decomposes into the direct sum
\begin{equation}
\label{Kasparek}
\P(n) = \bigoplus \colorop {\P}(c_1,\ldots,c_n;c),
\end{equation}
where the summation runs over all colors $c, \Rada c1n \in {\bf C}$. 
We require the decomposition~(\ref{Kasparek}) to be, in the
obvious sense, 
$\Sigma_n$-equivariant. More crucially, we demand the following.

Let $x \in \colorop{\P}(\Rada c1n;c)$ and $x_i\in
\colorop{\P}(\rada{d^i_1}{d^i_{k_i}};d_i)$, $1\leq i \leq n$. Then we
require that the non-triviality of the 
composition $x(\Rada x1n)$ implies that
\begin{equation}
\label{Chrochtatko}
d_i = c_i, \ \mbox {for $1\leq i \leq n$}, 
\end{equation}
in which case
\[
x(\Rada x1n) \in \colorop{\P}(\rada {d^1_1}{d^1_{k_1}},\ldots,
\rada {d^n_1}{d^n_{k_n}};c).
\]

The intuitive meaning of~(\ref{Chrochtatko}) is that one may plug
the element $x_i$ into the
$i$-th slot of the element $x$ if and only if the color of the output
of $x_i$ is the same as the color of the $i$-th input of $x$,
otherwise the composition is defined to be zero.

An example is provided by the {\em colored endomorphism operad\/} $\End_\bfU$
on a `colored' chain complex $\bfU = \bigoplus_{c\in \bfC} \bfU_c$
where we put, in~(\ref{Kasparek}), 
\[
\colorop{\End_\bfU}(c_1,\ldots,c_n;c) := {\rm Hom}(\orada
{U_{c_1}}{U_{c_n}},U_c). 
\]

{}From now on, by an operad we always mean
a \underline{colored} operad. 
Of course, any ordinary operad is
a ${\bf C}$-colored operad, with ${\bf C} = \{\rm Point\}$.
On the other hand, an arbitrary colored operad can be interpreted as an
ordinary operad defined over a certain semisimple algebra. To be more
precise, let $K := {\bfk_1} \oplus \cdots \oplus {\bfk_C}$ be the direct sum of copies
of the ground field indexed by elements of ${\bf C}$. It is easy
to see that ${\bf C}$-colored operads are the same as $K$-operads in
the sense of~\cite{ginzburg-kapranov:DMJ94}, compare 
also~\cite[Remark~1.97]{markl-shnider-stasheff:book}. 
This explains why most of the concepts of the theory of  
ordinary operads like 
free operads, ideals, presentations, etc., applies to the
colored case as well.

To set up the notation, we explain in detail what we mean by free
colored operads on colored collections.
A ${\bf C}$-colored collection is, roughly speaking, what remains if
one forgets composition operations of a colored operad, that is, a
system of ${\bf k}$-vector spaces
\[
E = \left\{
\colorop E(c_1,\ldots,c_n;c);\
c,\ c_1,\ldots,c_n \in \bfC,\ n \geq 1
\right\}
\]
together with a right action of the symmetric group $\Sigma_n$ on 
\[
E(n) := \bigoplus_{c,\ c_1,\ldots,c_n \in \bfC} \colorop E(c_1,\ldots,c_n;c)
\]
which permutes the colors in the obvious way. The {\em free operad functor\/}
$\Gamma(-)$ is then the left adjoint to the forgetful functor from the
category of colored operad to the category of these colored collections. 
The operad $\Gamma(E)$ is called the {\em free operad\/} on the colored
collection $E$.
If $E$ is linearly generated by a system of
generators, say $E = {\it Span}_{\bfk}(g_1,g_2,\ldots)$, we will
simplify the notation and write $\Gamma(g_1,g_2,\ldots)$
instead of $\Gamma({\it Span}_{\bfk}(g_1,g_2,\ldots))$.

\begin{example}
\label{certik_Pik}
{\rm\
The two-colored operad ${{{\Ass_\btb}}}$, with $\bfC := \{v,w\}$, 
describing diagrams
$(V,\mu) \stackrel f{\to} (W,\nu)$, where $(V,\mu)$ and $(W,\nu)$ are 
associative algebras and $f$
their homomorphism, has the presentation:
\begin{equation}
\label{Slune}
{{{\Ass_\btb}}} := 
\frac{\Gamma(\mu,\nu,f)}{(\mu(\mu \otimes \id) = \mu(\id \otimes \mu),\
\nu(\nu \otimes \id) = \nu(\id \otimes \nu),\ f\mu = \nu(f \otimes f))},
\end{equation}
where 
\[
\mu \in \colorop {{{{\Ass_\btb}}}}(v,v;v),\ \nu \in 
\colorop{{{{\Ass_\btb}}}}(w,w;w)
\mbox { and } f \in \colorop{{{\Ass_\btb}}}(v;w).
\]
In~(\ref{Slune}), $\Gamma(\mu,\nu,f)$ is 
the free operad on $\mu$ (a~generator for
the multiplication on $V$), $\nu$ (a~generator for the multiplication
on $W$) and $f$ (a generator for the map), modulo the ideal generated
by the axioms which have to be satisfied -- the associativity of $\mu$
and $\nu$ together with the axiom saying that $f$ is a homomorphism.

The diagram $(V,\mu) \stackrel f{\to} (W,\nu)$ 
as above is then the same as a
homomorphism $A: {{{\Ass_\btb}}} \to \End_{V,W}$ 
from the colored operad ${{{\Ass_\btb}}}$ to the colored operad
$\End_{V,W}:= \End_\bfU$ where $\bfU := U_v \oplus U_w$ with $U_v :=
V$ and $U_w := W$.
}
\end{example}

In the above example we in fact used the following

\begin{convention}
\label{koupil_jsem_si_GPS}
We will sometimes denote elements of operads  by the
same symbols as the corresponding multilinear maps of vector spaces.
\end{convention}

The category of operads admits coproducts. We call the coproduct of
operads $\calP_1$ and $\calP_2$ the {\em free product\/} and denote it
$\calP_1 * \calP_2$. 
If $\calP_1$ and $\calP_2$ are represented as
quotients of free operads, $\calP_s = \Gamma(X_s)/(R_s)$, $s=1,2$,
then $\calP_1 * \calP_2$ is represented by the quotient  $\Gamma(X_1,X_2)/(R_1,R_2)$,
where $(R_1,R_2)$ denotes the ideal
generated by $R_1 \cup R_2$.

\begin{example}
\label{jsem_zvedav_jestli_ve_Zbraslavicich_poletim}
{\rm\
This example generalizes Example~\ref{certik_Pik} to arbitrary algebras.
Given an (ordinary) operad $\P$, then the (colored) operad
\begin{equation}
\label{dve_prasatka}
\P_{\btb} := \frac{\P_{\bullet_1} * \P_{\bullet_2} *  \Gamma(f)}
                  {(fp_{\bullet_1} = p_{\bullet_2} f^{\ot n})},
\end{equation}
where $\P_{\bullet_1}$
(resp.~$\P_{\bullet_2}$) denotes a copy of $\P$ `concentrated' at the
first (resp.~the second) space and 
$fp_{\bullet_1} = p_{\bullet_2} f^{\ot n}$ means that, for each
$n\geq 1$ and each
operation $p \in \P(n)$, the arrow $f$ is a $p$-homomorphism, 
describes homomorphisms of $\P$-algebras. 
}
\end{example}

Example~\ref{jsem_zvedav_jestli_ve_Zbraslavicich_poletim} 
can be further generalized as follows. 
Let $S$ be a finite diagram, i.e.~a category with finite
number of objects and maps, and $\P$ an ordinary operad. The
finiteness assumption on the diagram
is not really necessary, but it will simplify
the exposition and we will always assume it.

Then there exist an operad $\P_{\uS}$ describing $S$-diagrams of
$\P$-algebras. The operad $\P_{\uS}$ can be constructed 
by taking generators for $\P$-algebra structures at each 
vertex of the diagram 
$S$, plus generators for the arrows of $S$, and then modding out
the ideal generated by the axioms of the $\P$-structures, 
the axioms saying that the arrows
of the diagram $S$ are $\P$-homomorphisms and the axioms 
expressing the commutativity relations inside $S$.

We need also the relative version of the above construction.
For a couple $(S,D)$ of a diagram $S$ and its subdiagram
$D$, by an  $(S,\uD)$-diagram of $\P$-algebras we mean an 
$S$-diagram of differential graded vector spaces with the property
that each vector space corresponding to a vertex of the
subdiagram $D$ is equipped with a $\P$-algebra structure such
that each arrow belonging to 
the subdiagram $D$ is a homomorphism of these
$\P$-algebras. 
Thus by underlying we indicate that the corresponding part of the
diagram carries a $\P$-structure.

\begin{definition}
\label{Krtecek}
Let  $\P_{(S,{\underline D})}$ denote the operad that describes
$(S,\uD)$-diagrams of $\P$-algebras.
The {\em canonical map\/} $\phi = \phi_\CD : \P_{(S,\uD)} \to
\P_{{\underline S}}$ is defined by forgetting the $\P$-structures on the
complement of $D$ in $S$.
\end{definition}

Notice that the map $\phi$ in the above definition goes in the
opposite direction than one would expect.
Observe that $\P_{(D,\uD)} = \P_{\uD}$, while $\P_{(D,\underline
\emptyset)} = {\bf 1}_{\uD}$,
the operad describing $D$-diagrams of graded vector spaces, i.e.~algebras
with no operations.

\begin{definition}
\label{Dva_kasparci_a_certik}
We say that the couple $(S,D)$ has
the (classical) {\em 
$\P$-extension property\/} if there exist a left inverse $\psi :
\P_{{\underline S}} \to \P_{(S,\underline D)}$ of the canonical map 
$\phi: \P_{(S,\underline D)} \to \P_{{\underline S}}$ of 
Definition~\ref{Krtecek},
\begin{equation}
\label{Usacek}
\psi \circ \phi = \id.
\end{equation}
\end{definition}

\begin{proposition}
\label{Zajicek_v_sve_jamce}
If $(S,D)$ has the $\P$-extension property, then each $(S,\uD)$-diagram
of $\P$-algebras can be extended to an $S$-diagram of $\P$-algebras.
\end{proposition}

\noindent 
{\bf Proof.}
We must show that, given a map $A_\CD: \P_{(S,\uD)} \to \End_\bfU$
representing an $(S,\uD)$-diagram of $\P$-algebras, there exists 
a map $A_\uS:\P_{(S,\uD)} \to \End_\bfU$ representing an $S$-diagram of
$\P$-algebras such that the diagram
\begin{center}
\VVtriangle{\P_{(S,\uD)}}{\P_\uS}{\End_\bfU}{\phi}{A_\CD}{A_\uS}
\end{center}
is commutative. The map $A_\uS:= A_\CD\circ \psi$ clearly has the
desired property.%
\qed

Let us formulate the following conjecture.

\begin{conjecture}
The implication of Proposition~\ref{Zajicek_v_sve_jamce} can be reversed,
i.e.~if any $(S,\uD)$ diagram of $\P$-algebras can be extended to an
$S$-diagram of $\P$-algebras, then the couple $(S,\uD)$ has the 
$\P$-extension property. 
\end{conjecture}

\begin{example}[related to {\bf (M1)}]
\label{7}
{\rm
Let $S$ be the diagram
\begin{center}
{
\unitlength=.45pt
\begin{picture}(180.00,40.00)(80.00,30.00)
\put(20.00,10.00){\vector(-3,1){1.00}}
\put(80.00,52.00){$f$}
\put(80.00,-22.00){$g$}
\put(160.00,30.00){\vector(3,-1){1.00}}
\bezier{100}(20.00,10.00)(90.00,-10.00)(160.00,10.00)
\bezier{100}(20.00,30.00)(90.00,50.00)(160.00,30.00)
\put(180.00,20.00){\makebox(0.00,0.00){$\circ$}}
\put(0.00,20.00){\makebox(0.00,0.00){$\bullet$}}
\put(187.00,20.00){\makebox(0.00,0.00)[l]{, $fg=1$ and $gf =1$,}}
\end{picture}}
\end{center}
\vskip .8cm

\noindent 
and let the subdiagram $D$ consist of the `solid vertex' and no morphisms.
Then 
\begin{equation}
\label{Hopsalek}
\P_{(S,\underline D)} := \frac{\P_{\bullet} * \Gamma(f,g)}{(fg=1,\ gf=1)},
\end{equation}
where $\P_{\bullet}$ denotes the copy 
of $\P$ `concentrated' at the vertex
$\bullet$. In a similar manner,
\[
\P_{\underline S}:=\frac{ \P_{\bullet}* \P_{\circ} * \Gamma(f,g)}{(
fg=1,\ gf=1,\ fp_{\bullet} = p_{\circ} f^{\ot n},\ gp_{\circ} =
p_{\bullet } g^{\ot n})},
\]
where we use notation as in~(\ref{dve_prasatka}).
The canonical
map $\phi : \P_{(S,\underline D)} \to \P_{\underline S}$ is given
by $\phi(f): = f$, $\phi(g):= g$ and $\phi(p_{\bullet}):=
p_{\bullet}$. 
Let us consider a map $\psi: \P_{\underline S} \to \P_{(S,\underline D)}$
defined by
\begin{equation}
\label{Pik}
\psi(f):= f,\
\psi(g):= g,\
\psi(p_{\bullet}):= p_{\bullet}
\mbox { and }
\psi(p_{\circ}):= f p_{\bullet} g^{\ot n}.
\end{equation}
We must verify that this definition is correct. Clearly, $\psi(fg)
= fg$  and $\psi(gf)= gf$, thus $\psi$ maps the axioms $fg=1$ and $gf
= 1$ to themselves.
Next, we have $\psi(f p_\bullet) = f p_\bullet$ and
\[
\psi(p_\circ f^{\ot n})= f p_\bullet g^{\ot n}  f^{\ot n} = 
f p_\bullet {(gf)}^{\ot n}=
\mbox { (because $gf = 1$ in $\P_{(S,\uD)}$) } = f p_\bullet,
\]
thus $\psi$ respects 
the axiom $f p_\bullet = p_\circ f^{\ot n}$. In a similar
manner we verify that $\psi(g p_\circ) = p_\bullet g^{\ot n} =
\psi(p_\bullet g^{\ot n})$, which finishes the verification that
$\psi$ is well-defined.

It is immediately clear by checking on the generators
that $\psi \circ \phi = \id$, so the couple
$(S,D)$ has the $\P$-extension property. It is very easy to
verify that also $\phi\circ \psi = \id$, thus the map $\phi$ is in
fact an isomorphism. 
}
\end{example}

\begin{example}[related to (M3)]
\label{9}
{\rm
Let $S$ be the diagram
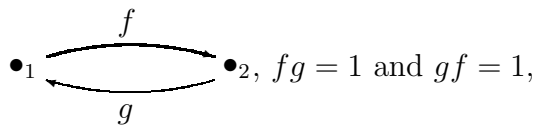
\begin{center}
{
\unitlength=.45pt
\begin{picture}(180.00,40.00)(70.00,30.00)
\put(20.00,10.00){\vector(-3,1){1.00}}
\put(80.00,52.00){$f$}
\put(80.00,-22.00){$g$}
\put(160.00,30.00){\vector(3,-1){1.00}}
\bezier{100}(20.00,10.00)(90.00,-10.00)(160.00,10.00)
\thicklines
\bezier{100}(20.00,30.00)(90.00,50.00)(160.00,30.00)
\put(180.00,20.00){\makebox(0.00,0.00){$\bullet_2$}}
\put(193.00,20.00){\makebox(0.00,0.00)[l]{, $fg=1$ and $gf=1$,}}
\put(0.00,20.00){\makebox(0.00,0.00){$\bullet_1$}}
\end{picture}}
\end{center}
\vskip 1cm

\noindent 
and let $D$ be the `thick' subdiagram. Then
\begin{equation}
\label{meda_Pusik}
P_{(S,\underline D)}: = 
\frac{ \P_{\bullet_1}* \P_{\bullet_2} * \Gamma(f,g)}{(fg=1,\ gf=1,\
fp_{\bullet_1} = p_{\bullet_2} f^{\ot n})}
\end{equation}
and
\[
P_{\underline S}: = 
\frac{ \P_{\bullet_1}* \P_{\bullet_2} * \Gamma(f,g)}{(fg=1,\ gf=1,\
fp_{\bullet_1} = p_{\bullet_2} f^{\ot n},\
gp_{\bullet_2} =
p_{\bullet_1 } g^{\ot n})}.
\]
The canonical map $\phi: \P_\CD \to \P_\uS$, given by 
$\phi(f): = f$, $\phi(g):= g$, $\phi(p_{\bullet_1}):=
p_{\bullet_1}$ and $\phi(p_{\bullet_2}):=
p_{\bullet_2}$ is again immediately seen to be  an isomorphism.
}
\end{example}

We are going to give an example in which  the map $\phi$ has a left
inverse $\psi$, but $\phi$ is not an isomorphism. For this example we
need to introduce augmented operads. This notion will be useful also later in
the paper.
Let ${\bf 1}$ be the trivial (initial) ${\bf C}$-colored operad defined by
\begin{equation}
\label{jsem_zvedav_jestli_poletim_s_Vosou}
\colorop {\bf 1}(c_1,\ldots,c_n;c) := 
\cases{\bfk}{for $n=1$, $c_1=c$, and}{0}{otherwise.}
\end{equation}
For a ${\bf C}$-colored operad $\calP$, let $\eta : {\bf 1} \to \calP$
be the unique map induced by the unit of $\calP$.

\begin{definition}
A colored operad $\calP$ is {\em augmented\/} if there exists a
morphism $\epsilon : \calP \to {\bf 1}$ such that $\epsilon\eta =
{\it id}_{\bf 1}$. The ideal $\overline \calP := \Ker(\epsilon)$ is called
the {\em augmentation ideal\/}.
\end{definition}

Free (colored) operads are augmented.  Ordinary quadratic
operads~\cite{ginzburg-kapranov:DMJ94} (such as $\Ass$, $\Comm$ and
$\Lie$) are also augmented. 
If $\calP$ is an ordinary operad, $\calP_\CD$ is augmented
if and only if $\calP$ is augmented and the diagram
$S$ has no loops, that is, no nontrivial identity of arrows 
of the form $F = \id$
(such as $fg =1$) is required in $S$. 
Therefore, colored operads considered in
Examples~\ref{9}, \ref{Salvatore} 
and~\ref{tyden_pred_mistrovstvim} are not augmented.

\begin{example}[due to P.~Salvatore]
\label{Salvatore}
{\rm
Let $S$ be the diagram in which
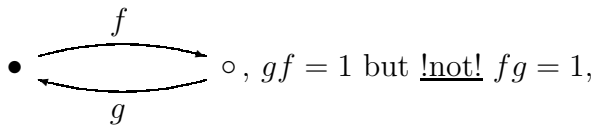
\begin{center}
{
\unitlength=.45pt
\begin{picture}(180.00,40.00)(90.00,30.00)
\put(20.00,10.00){\vector(-3,1){1.00}}
\put(80.00,52.00){$f$}
\put(80.00,-22.00){$g$}
\put(160.00,30.00){\vector(3,-1){1.00}}
\bezier{100}(20.00,10.00)(90.00,-10.00)(160.00,10.00)
\bezier{100}(20.00,30.00)(90.00,50.00)(160.00,30.00)
\put(180.00,20.00){\makebox(0.00,0.00){$\circ$}}
\put(191.00,20.00){\makebox(0.00,0.00)[l]{, $gf=1$ but
\underline{!not!}
$fg=1$,}}
\put(0.00,20.00){\makebox(0.00,0.00){$\bullet$}}
\end{picture}}
\end{center}
\vskip 1cm

\noindent 
and let $D$ be the  solid vertex. Then 
\[
\P_{(S,\uD)} := \frac{\P_\bullet * \Gamma(f,g)}{(gf=1)} \hskip .5cm
\mbox { and }\hskip .5cm
\P_\uS := \frac{\P_\bullet* \P_\circ * \Gamma(f,g)}{(gf=1,\ fp_\bullet =
p_\circ f^{\ot n},\ gp_\circ = p_\bullet g^{\ot n})}.
\]

The canonical map $\phi: \P_\CD \to \P_\uS$ is given by 
$\phi(p_\bullet):= p_\bullet$, $\phi(f)=
f$ and $\phi(g)= g$. Let us try to define the map $\psi : \P_\uS \to
\P_{(S,\uD)}$ by $\psi(f):= f$, $\psi(g):= g$, $\psi(p_{\bullet}):=
p_{\bullet}$ and
\begin{equation}
\label{tak_malo_casu}
\psi(p_{\circ}):= f p_{\bullet} g^{\ot n}
\end{equation}
(i.e~by formally the same formula as~(\ref{Pik})).
Unfortunately, the map $\psi$ is not well-defined. 
Equation~(\ref{tak_malo_casu}) gives, for the unit $1 \in \calP(1)$,
$\psi(1_\circ) = f 1_\bullet g = fg$, therefore $\psi(1_\circ) \not=
1_\circ$ (recall we did not assume $fg=1_\circ$), thus $\psi$ is not an operad
morphism.

For {\em augmented operads\/}, this subtlety can be fixed  
by requiring~(\ref{tak_malo_casu}) only for $p \in \overline {\calP}$.
It is then clear that $\psi\circ \phi = \id$. 
On the other hand, the map $\phi$ is not  an isomorphism, because the
equivalence class of $p_\circ$ in $\P_\uS$ is clearly not in the image
of $\phi$.
}
\end{example}

The `move' corresponding to the diagram of Example~\ref{Salvatore} is:

\begin{itemize}
\item[{\bf (S)}]
Let $\calP$ be augmented. Then, for any strongly 
homotopy $\P$-algebra ${\bf V}$, chain complex $W =
(W,\pa)$ and a chain map $f : (V,\pa) \to (W,\pa)$ that has a left
chain-homotopy inverse,  there exist a sh{}
$\P$-structure $\bfW$ on the chain complex
$(W,\pa)$ and a strongly homotopy $\P$-structure ${\bf f}:\bfV \to \bfW$
on the map $f$.
\end{itemize}

Let us emphasize again that the augmentation of
$\calP$ is crucial here. This excludes  
{\em algebras with units\/}, because operads describing these algebras (such
as $\calP = \mbox{\it U \hskip -1.5mm \Ass}$ 
for associative algebras with unit)
are never augmented.

\begin{example}[{\bf related to (M4)}]
{\rm
We leave it to the reader to analyze the situation when $S$ is the
diagram
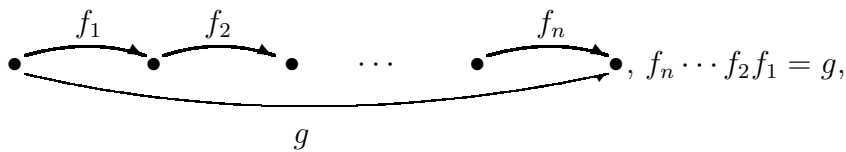
\begin{center}\hskip 1cm
{
\unitlength=.35pt
\begin{picture}(650.00,80.00)(200.00,30.00)
\put(620.00,25.00){\vector(2,1){20.00}}
\bezier{700}(10.00,35.00)(310.00,-35.00)(640.00,35.00)
\put(310.00,-35.00){\makebox(0,0){$g$}}
\thicklines
\put(120.00,65.00){\vector(2,-1){20.00}}
\bezier{200}(10.00,55.00)(80.00,75.00)(140.00,55.00)
\put(80.00,80.00){\makebox(0,0)[b]{$f_1$}}
\put(260.00,65.00){\vector(2,-1){20.00}}
\bezier{200}(160.00,55.00)(220.00,75.00)(280.00,55.00)
\put(220.00,80.00){\makebox(0,0)[b]{$f_2$}}
\put(620.00,65.00){\vector(2,-1){20.00}}
\bezier{200}(510.00,55.00)(580.00,75.00)(640.00,55.00)
\put(580.00,80.00){\makebox(0,0)[b]{$f_n$}}
\put(650.00,45.00){\makebox(0.00,0.00){$\bullet$}}
\put(661.00,45.00){\makebox(0.00,0.00)[l]{, $f_n \cdots f_2f_1 =g$,}}
\put(390.00,45.00){\makebox(0.00,0.00){$\cdots$}}
\put(500.00,45.00){\makebox(0.00,0.00){$\bullet$}}
\put(300.00,45.00){\makebox(0.00,0.00){$\bullet$}}
\put(150.00,45.00){\makebox(0.00,0.00){$\bullet$}}
\put(0.00,45.00){\makebox(0.00,0.00){$\bullet$}}
\end{picture}}
\end{center}

\vskip 9mm
\noindent 
and $S$ is the `solid' subdiagram, that is, to describe operads
$\P_\CD$, $\P_\uS$ and prove that the canonical map
$\phi : \P_\CD \to \P_\uS$ is an isomorphism.
}
\end{example}

\begin{example}
{\rm
Let us give an example which {\em does not\/} have the extension property. Let
$S$ be the diagram
\begin{center} 
\hskip 1cm
{
\unitlength=.45pt
\begin{picture}(180.00,40.00)(80.00,30.00)
\put(80.00,52.00){$f$}
\put(180.00,20.00){\makebox(0.00,0.00){$\circ$}}
\put(0.00,20.00){\makebox(0.00,0.00){$\bullet$}}
\thinlines
\put(20.00,20.00){\vector(1,0){140}}
\end{picture}}
\end{center}
\vskip .8cm

\noindent  
and let $D$ be the subdiagram consisting of the solid vertex
$\bullet$. Then
\[
\P_{(S,\underline D)}:= {\P_{\bullet} * \Gamma(f)},\
\P_{\underline S}:= \frac{\P_{\bullet} * \P_{\circ} * 
\Gamma(f)}{(fp_{\bullet} = p_{\circ} f^{\ot n})},
\]
and the canonical map $\phi: \P_\CD \to \P_\uS$ 
is given by $\phi(p_{\bullet}):=
p_\bullet$ and $\phi(f):= f$. 

Suppose that $\psi \circ \phi
= \id$ for some map $\psi : \P_\uS \to \P_\CD$. Then clearly 
$\psi(p_\bullet) = p_\bullet$ and $\psi(f)=f$. The element $\omega :=
\psi(p_\circ)$ must satisfy $f p_\bullet = \omega f^{\ot n}$ in
$\P_\CD$, otherwise $\psi$ will not be well-defined. But it is
immediate to see that there is no such $\omega$, thus $\phi$ has no
left inverse.

On the algebra level this means the existence of a $\P$-algebra
$V$, a vector space $W$ and a linear map $f : V \to W$ 
for which there is no
$\P$-algebra structure on $W$ such that $f$ is a homomorphism.

To construct an example, take any $\P$-algebra $V$ with a
subspace $I \subset V$ that is {\em not\/} an ideal. Put $W := V/I$
(quotient of vector spaces). Then the canonical
projection $f: V \to W$ clearly has the above property.
}
\end{example}

Let us close this  section with the following problem.

\noindent
{\bf Problem.}
{\em
Give a characterization of pairs $(S,D)$ that have the $\P$-extension
property. Is there a (finite) list of `moves' such that $(S,D)$ has
the $\P$-extension property if and only if $S$ is obtained from $D$ by
consecutively applying moves from the list? 

To which extent the
$\calP$-extension property depends on $\calP$? Is it, for example,
true that if the pair $(S,D)$ has the $\calP$-extension property for
some nontrivial augmented operad $\calP$, then it has the
$\calP$-extension property for an arbitrary augmented $\calP$?
}

In all examples with the $\P$-extension property we know, the bigger
diagram $S$ is obtained from $D$ by the moves listed above.

\section{Resolutions, homotopy diagrams and extensions}
\label{kam_nas_vystehuji?}

As we mentioned in the introduction, by a strongly homotopy $\P$-algebra we
mean an algebra over a minimal, or at least cofibrant, 
resolution of the operad
$\P$. The same approach translates to the word of colored
operads. Let us illustrate this on the following example.

\begin{example}
\label{tenhle_clanek_je_nekonecna_story}
{\rm
Let ${{\Ass_\btb}}$ be the operad from Example~\ref{certik_Pik} in
presentation~(\ref{Slune}). The minimal
model of ${\Ass_\btb}$ is given by (see~\cite{markl:ho})
\[
{\Ass_\btb} \stackrel{\alpha}{\longleftarrow}
\left(
\Gamma\lbigbrace
\mu_2,\mu_3,\mu_4,\ldots,f_1,f_2,f_3,\ldots,
\nu_2,\nu_3,\nu_4,\ldots \rbigbrace,\pa
\right),
\]
where $\alpha(\mu_2) = \mu$, $\alpha(\nu_2)= \nu$, $\alpha(f_1)= f$,
while $\alpha$ is trivial on the remaining generators.
The differential is given by the formulas
\begin{eqnarray*}
\pa (\mu_n)
&:=&
\sum_{\doubless{i+j=n+1}{i,j \geq 2}} \sum_{s=0}^{n-j}
(-1)^{j+s(j+1)} \mu_{i}(\id^{\ot s} \ot \mu_j \ot
\id^{\ot i-s-1}),
\\
\pa (f_n) &:=&
\sum_{k=2}^n \sum_{\prada{r_1}{r_k} =n} 
(-1)^{\sum_{1\leq i< j \leq k}r_i(r_j+1)}
\nu_k(\orada{f_{r_1}}{f_{r_k}})-
\\
&&\hskip 1cm- \hskip -3mm
\sum_{\doubless{i+j=n+1}{i,j \geq 2}}
\sum_{s=0}^{n-j}
(-1)^{n+s(j+1)}
f_{i}(\id^{\ot s} \ot \mu_j \ot
\id^{\ot i-s-1}),
\\
\pa (\nu_n)
&:=&
\sum_{\doubless{i+j=n+1}{i,j \geq 2}} \sum_{s=0}^{n-j}
(-1)^{j+s(j+1)} \nu_{i}(\id^{\ot s} \ot \nu_j \ot
\id^{\ot i-s-1}),
\end{eqnarray*}
Comparing this with formulas~(\ref{An}) and~(\ref{Fn}),
we conclude that a strongly homotopy ${{\Ass_\btb}}$-algebra on a
colored chain complex $V\oplus W$ is given by an $\ainfty$-structure on
$V$ (generators $\mu_2,\mu_3,\ldots$), an $\ainfty$-structure on
$W$ (generators $\nu_2,\nu_3,\ldots$), and a sh{} morphism of these two
structures (generators $f_1,f_2,\ldots$).
}
\end{example}

The above example motivates the following definition.

\begin{definition}
\label{nikdy_se_toho_nezbavim}
A {\em strongly homotopy $(S,\uD)$-diagram of $\P$-algebras\/} 
is an algebra over a
cofibrant resolution of the operad $\P_{(S,\underline D)}$.
\end{definition}

At this stage, it is necessary to say more
precisely what we mean by a cofibration.

\begin{definition}
\label{cof}
By an {\em elemental cofibration\/} we mean a map of dg-operads of 
the form $\iota : (\S,\pa) \hookrightarrow (\S * \Gamma(E),\pa)$ 
such that $E$ decomposes as $E = \bigoplus_{n \geq 0}E_n$ and
\begin{equation}
\label{Opicak_Fuk}
\pa(E_n) \subset S* \Gamma(E)_{< n},
\mbox { for each $n \geq 0$,}
\end{equation}
where $\Gamma(E)_{< n}$ denotes the suboperad of $\Gamma(E)$ generated by
$\bigoplus_{j < n} E_j$.

An operad $\S$ is {\em elementally cofibrant} 
if the canonical map $\eta : {\bf
1} \to \calS$ is an elemental cofibration or, equivalently, if $\calS$ is
of the form $(\Gamma(E),\pa)$, where $E$ decomposes as above and
\begin{equation}
\label{Strihatko}
\pa(E_n) \subset \Gamma(E)_{< n},
\mbox { for each $n \geq 0$}.
\end{equation}
\end{definition}

The definition above is an `operadic' analog of the Koszul-Sullivan
algebra in rational homotopy theory, see~\cite{halperin:lect}.
We believe that
there exists a closed model category (CMC) structure on
the category of (colored) operads such that fibrations are epimorphisms, weak
equivalences are homology isomorphisms (quisms) and cofibrations are 
retractions of elemental cofibrations.
We are, however, not going to discuss this CMC structure here, because
it will be enough for our purposes to prove a certain factorization
property (Theorem~\ref{factorization}) and a~lifting lemma
(Lemma~\ref{ll}) for cofibrations. Let us make first the following

\begin{convention}
By a cofibration we will 
\underline{alwa}y\underline{s} \underline{mean} an elemental one.
\end{convention}

The factorization property says:

\begin{theorem}
\label{factorization}
For an arbitrary map  $f: \S  \to {\cal Q}$ of differential colored operads,
there exists a cofibration $\iota : \S \hookrightarrow \calR$ and a homology
isomorphism $\alpha: \calR \to {\cal Q}$ making the diagram
\begin{center}
\Vtriangle{\S}{{\cal Q}}{\calR}f{\iota}{\alpha}
\end{center}
commutative.
\end{theorem}

We will need also the following `lifting lemma:'

\begin{lemma}
\label{ll}
Consider the following commutative diagram of solid arrows between
colored operads:
\[
\squarewithdotedmensi{\calS}{\calX}{\calQ}{\calY}f{\iota}pgh
\]
in which $p$ is an epimorphism and homology isomorphism, and $\iota$
an (elemental) cofibration. Then there exists a morphism $h: \calQ
\to \calX$ such that $h\iota = f$ and $ph = g$.
\end{lemma}

Theorem~\ref{factorization} and Lemma~\ref{ll} are proved in the Appendix. 
Theorem~\ref{factorization} implies that each ${\bf C}$-colored
operad $\calP$ indeed has a {\em cofibrant model\/} -- a cofibrant
operad $\calR$ together with a homology isomorphism $\alpha : \calR \to
\calP$. To see this, apply the theorem to
$\calQ = \calP$, $\calS = {\bf 1}$ and $f = \eta$, 
where ${\bf 1}$ is the trivial (initial) ${\bf
C}$-colored operad defined in~(\ref{jsem_zvedav_jestli_poletim_s_Vosou})
and $\eta : {\bf 1} \to \calP$ is induced by the unit of $\calP$.

This provides a resolution for an arbitrary
colored operad, but gives no control of the
size of this resolution. We show below that some more special
colored operads have very small resolutions.

\begin{definition}
\label{zitra_jedu_do_Oberwolfachu}
A {\em bigraded minimal model\/} of a colored operad $\calP$ 
with trivial differential is given by the following data:
a bigraded collection $Z_*^*=Z^0_*\op Z_*^1\op\cdots$,
a degree -1 differential $\pa$ on $ \Gamma (Z)$ and
a map $\rho : (\Gamma(Z),d)\to (\calP,0)$ of differential operads
such that the following conditions are satisfied:
\begin{enumerate}
\item[(i)]
$\pa$ is {\em minimal\/} in the sense that $\pa(Z)$ consist of decomposable
elements of $\Gamma(Z)$,
\item[(ii)]
$\pa$ decreases the upper grading of $\Gamma(Z)$ by one, $\pa(Z^i) \subset
\Gamma(Z)^{i-1}$,
\item[(iii)]
$\rho|_{Z^{\geq 1}}=0$ and $\rho$ induces an isomorphism
$H^0(\rho):H^0(\Gamma(Z),\pa) \cong \calP$, and
\item[(iv)]
$H^{\geq 1}(\Gamma(Z),d)=0$.
\end{enumerate}
We call the grading induced by the ``upper'' grading of $Z$ the
{\em TJ-grading\/} (from Tate-Josefiak). 
\end{definition}

While we proved in~\cite{markl:zebrulka} that any ordinary
(non-colored) operad $\calP$ with $\calP(0) =0$ and
$\calP(1) \cong \bfk$ admits a bigraded minimal model, 
we do not know any reasonably large class of
colored operads with this property.

It immediately follows from (i) and (iii) of 
Definition~\ref{zitra_jedu_do_Oberwolfachu}
that any operad having a bigraded minimal model must be augmented.
Thus for example the operad $\Iso$ introduced later
in Example~\ref{tyden_pred_mistrovstvim}
does not admit a minimal bigraded model. On the other hand, consider the
diagram $D$ consisting of one vertex and three endomorphisms $f$, $g$, $h$ such
that $fg = h$, $gh = f$ and $hf{} = g$. Then $\calP_D$ does not admit
a bigraded minimal model although  $\calP_D$ is augmented for each augmented
$\calP$. Fortunately, for the scope of applications of this paper, the
following statement is sufficient.

\begin{proposition}
\label{TJ}
Let $\calP$ be an ordinary operad with $\calP(0) = 0$, $\calP(1)
\cong \bfk$ and trivial differential. Then
both $\calP$ and the colored operad 
$\calP_{\btb}$ admit a minimal bigraded model. Moreover, this
model is unique up to an isomorphism
and each cofibrant resolution is isomorphic to the 
free product of this bigraded model and an acyclic free operad.
\end{proposition}

\noindent 
{\bf Proof.}
The bigraded minimal model for $\calP$ was constructed
in~\cite[Theorem~2.5]{markl:zebrulka}. In~\cite[Theorem~7]{markl:ho}
we constructed the minimal model for $\calP_\btb$ for an arbitrary
dg-operad with $\calP(0) = 0$ and $\calP(1) \cong \bfk$.
It is immediate to see that this model was 
in fact minimal bigraded when $\calP$ had trivial differential.%

The second part of the proposition follows from a word-by-word
translation of the proof of Theorem~2.2 
in~\cite[page~282]{sullivan:Publ.IHES77}.
\qed

For the proof of Proposition~\ref{zitra_brzo_vstavam} formulated below,
we need the following form of K\"unneth formula.

\begin{theorem}
\label{easy_kunneth}
Given augmented dg-operads $\calS$ and $\calQ$, 
one has a functorial isomorphism
\[
H_*(\calS)* H_*(\calQ) \cong H_*(\calS * \calQ).
\]
\end{theorem}

The proof of Theorem~\ref{easy_kunneth} is postponed to the Appendix.
The following proposition says that each cofibrant model of the operad
$\calP$ or $\calP_\btb$ can be equipped with a TJ-grading. We will
need this fact in the proof of Theorem~\ref{takhle_prosvihnout_seky}.

\begin{proposition}
\label{zitra_brzo_vstavam}
Let $\calP$ be as in Proposition~\ref{TJ} and let $(\Gamma(E),\pa)$ be
a cofibrant model of $\calP$ or $\calP_\btb$. Then  $(\Gamma(E),\pa)$
is isomorphic to a model of the form $(\Gamma(X),\pa)$ which satisfy
all conditions of Definition~\ref{zitra_jedu_do_Oberwolfachu} except
possibly the minimality~(i). We call $(\Gamma(X),\pa)$ the {\em
bigraded model\/}.
\end{proposition}
 
\noindent
{\bf Proof.} 
By Proposition~\ref{TJ}, $(\Gamma(E),\pa)$ is isomorphic to
\begin{equation}
\label{boli_mne_siska}
(\Gamma(Z),\pa)* \Gamma(A,\pa A) \cong \Gamma(Z \oplus A \oplus \pa A,\pa),
\end{equation}
where $\rho : (\Gamma(Z),\pa) \to (\calP,\pa=0)$ is the
bigraded minimal model. Put $X := Z \oplus A \oplus \pa A$ with
TJ-grading induced by TJ-grading of $Z$ and postulating $\pa
A$ (resp.~$A$) to be concentrated in TJ-grading $0$
(resp.~$1$). Condition~(ii) of
Definition~\ref{zitra_jedu_do_Oberwolfachu} is clearly satisfied. The
acyclicity of the operad $\Gamma(A,\pa A)$ is obvious -- 
there exists a very explicit contracting
homotopy for this operad. Conditions~(iii) and~(iv) can be verified by applying
Theorem~\ref{easy_kunneth} to the product~(\ref{boli_mne_siska}) of
augmented dg-operads.%
\qed
 
\begin{example}
{\rm
The resolution of the operad $\Ass_\btb$ described in
Example~\ref{tenhle_clanek_je_nekonecna_story} is in fact the bigraded
minimal model of $\Ass_\btb$, with $Z^n =
\span(\mu_{n+2},f_{n+1},\nu_{n+2})$, $n \geq 0$. 
In this example, the TJ grading coincides with the 
homological grading, which is always the case when $\calP$ and
therefore also $\calP_\btb$ is concentrated in degree $0$.
}
\end{example}

Consider a couple $(S,D)$ of a diagram and its subdiagram, 
and a cofibrant resolution $\alpha_{\CD}:
\R_{\CD} \to \P_{\CD}$ of the operad $\P_{\CD}$. By
Theorem~\ref{factorization} (with $\calS = \calR_\CD$, $\calQ =
\calP_{{\underline S}}$ and $f = \phi \circ \alpha_\CD$) 
there exist a cofibrant resolution 
$\alpha_{\underline S} : \R_{\uS} \to \P_{\uS}$ and a cofibration
$\iota:\R_{\CD} \to  \R_{\uS}$ such that the diagram
\begin{center}
\square{\R_{\CD}}{\R_{\uS}}{\P_{\CD}}{\P_{\uS}}%
        {\iota}{\alpha_{\CD}}{\alpha_{\uS}}{\phi}
\end{center}
commutes.

\begin{definition}
The couple $(S,D)$ has the {\em homotopy $\P$-extension property\/} if
there exist a left inverse $p: \R_{\uS} \to \R_{\SD}$ of the map
$\iota$ from the diagram above.
\end{definition}

Exactly as in the proof of Proposition~\ref{Zajicek_v_sve_jamce},
we show that the
$\P$-homotopy extension property implies that each strongly
homotopy $(S,\uD)$-diagram of $\P$-algebras can be extended to a
strongly homotopy $S$-diagram of $\P$-algebras.
The main theorem of the paper then reads.

\begin{theorem}
\label{main}
The couple $(S,D)$ has the homotopy extension property if and only if it
has the (classical) extension property of 
Definition~\ref{Dva_kasparci_a_certik}.
\end{theorem}

\noindent 
{\bf Proof.}
Suppose the couple $(S,D)$ has the classical extension property, that
is, that there exists a left inverse $\psi : \P_\uS \to \P_\CD$ of
the canonical map $\phi: \P_\CD \to \P_\uS$.
Consider the commutative diagram
\begin{center}
\squarewithdoted{\R_\CD}{\R_\CD}{\R_\uS}{\P_\CD}%
       {\id}{\iota}{\alpha_\CD}{\psi\alpha_\uS}{p} 
\end{center}
Since $\iota$ is a cofibration and $\alpha_\CD$ an epimorphism and 
homology isomorphism, by~Lemma~\ref{ll}
there exists a map $p : \R_\uS \to \R_\CD$ such
that $p\circ \iota = \id$, that is, $p$ is a left inverse of the map~$\iota$.

The opposite implication is clear. If $p \circ\iota = \id$ for some map
$p$, then 
\[
H(\alpha_\CD) \circ H(p) \circ H(\alpha_\uS)^{-1} : \P_\uS
\to \P_\CD
\] 
(the composition of the homology of maps and their inverses) 
is a left inverse of the
canonical map $\phi$.%
\qed

We close this section by a definition of the underlying
map of a sh{} morphism. As explained in 
Definition~\ref{nikdy_se_toho_nezbavim}, a sh{} $\calP$-homomorphism
of two sh{} $\calP$-algebras is an operadic morphism
$\bff :A_\btb : \calR_\btb \to \End_{U,V}$, where $\alpha_\btb : \calR_\btb
\to \calP_\btb$ is a cofibrant resolution of the bi-colored operad
$\calP_\btb$ introduced in~(\ref{dve_prasatka}). 

Let $\Gamma(u)$ be the free operad generated by an 
arrow $u : \bullet \to \bullet$. There certainly exists a (non-unique) 
operadic map
\begin{equation}
\label{dnes_jsem_blue}
\#: \Gamma(u) \to \calR_\btb
\end{equation}
such that $\alpha_\btb(\#(u)) = f$, where $f$ is the generator of
$\calP_\btb$ as in~(\ref{dve_prasatka}). Observe that, replacing
$\calR_\btb$ by an isomorphic resolution if necessary, we can always
make $\#$ an (elemental) cofibration.

\begin{definition}
\label{blue_mood}
The {\em underlying map\/} of the sh{} $\calP$-homomorphism $\bff
:A_\btb : \calR_\btb \to \End_{U,V}$ is the dg map $\#\bff := A_\btb
(\#(u)): (U,\pa) \to (V,\pa)$.

Conversely, a {\em strongly homotopy structure\/} on a dg homomorphism
$g : (U,\pa) \to (V,\pa)$ is given by a sh{} $\calP$-homomorphism
$\bff$ as above such that $g = \#\bff$.
\end{definition} 

The above definition involves a choice of the map $\#$, 
but as we argued in the introduction, in concrete situations
there is always a preferred choice. For instance, 
in Example~\ref{tenhle_clanek_je_nekonecna_story} one
may take $\#(u) := f_1$, thus the underlying map of an
$A_\infty$-morphism $(f_1,f_2,f_3,\ldots)$ is the dg map $f_1$ --
cf.~our Convention~\ref{koupil_jsem_si_GPS}.

\section{Structure of strongly homotopy diagrams}
\label{zase_starosti}

A little care is needed to apply the results of the previous section
to concrete situations. The difficulties can be summarized by saying that
strongly homotopy diagrams are often very complicated objects.
We illustrate this in Example~\ref{tyden_pred_mistrovstvim} on the
resolution $\Riso$ of the colored operad $\Iso$ describing
isomorphisms of chain complexes. The colored operad $\Riso$ will be
needed in Section~\ref{jeste_jeden_medvidek}. The main result of this
section is Lemma~\ref{sroubovacky}.

\begin{example}
\label{tyden_pred_mistrovstvim}
{\rm
Consider the diagram $S$ given by the following picture (we already
saw this diagram in Example~\ref{7}):
\begin{equation}
\label{hvezdicka}
{
\unitlength=.45pt
\begin{picture}(180.00,40.00)(80.00,30.00)
\put(20.00,10.00){\vector(-3,1){1.00}}
\put(80.00,52.00){$f$}
\put(80.00,-22.00){$g$}
\put(160.00,30.00){\vector(3,-1){1.00}}
\bezier{100}(20.00,10.00)(90.00,-10.00)(160.00,10.00)
\bezier{100}(20.00,30.00)(90.00,50.00)(160.00,30.00)
\put(180.00,20.00){\makebox(0.00,0.00){$\circ$}}
\put(0.00,20.00){\makebox(0.00,0.00){$\circ$}}
\put(187.00,20.00){\makebox(0.00,0.00)[l]{, $fg=1$ and $gf =1$.}}
\end{picture}}
\end{equation}
\vskip .5cm

\noindent 
An (ordinary) homotopy $S$-diagram is, by definition,  
an $S$-diagram in the homotopy category of chain complexes. 
For $S$ above, it is simply a chain homotopy equivalence,
that is, it consist of two arrows, $f$ and $g$, and
homotopies $h : gf \sim \id$ and $l : fg \sim \id$. This could be
algebraically expressed as a `resolution'
\[
\P_{(S,\underline \emptyset)} := \frac{\Gamma(f,g)}{(fg = 1,\ gf = 1)}
\stackrel {\alpha}{\longleftarrow}
(\Gamma(f,g,h,l),\partial),
\]
with $\alpha(f): = f$, $\alpha(g) := g$, $\alpha(h) = \alpha(l) =0$,
and the differential given by
\[
\pa(f) = \pa(g) := 0,\
\pa(h) := gf - 1
\mbox { and }
\pa(l) := fg -1.
\]
Unfortunately, the map $\alpha$ is not a homology
isomorphism. For example, the element $gl - hg \in\Gamma(f,g,h,l)$ is a
cocycle, mapped to zero by $\alpha$, that is not a coboundary. A
honest cofibrant model of $\P_{(S,\underline \emptyset)}$ 
can be described as follows:
\begin{equation}
\label{honest}
\P_{(S,\underline \emptyset)} \stackrel {\alpha}{\longleftarrow}
(\Gamma(f,g,h,l,f_2,g_2,f_3,g_3,\ldots),\paiso),
\end{equation}
with the differential $\pa = \paiso$ given by
\def\miniminus{{\hskip -.3mm - \hskip -.3mm}}
\def\miniplus{{\hskip -.3mm + \hskip -.3mm}}
\def\minidef{{\hskip -.2mm :=\hskip -.2mm}}
\[
\def\arraystretch{1.4}
\begin{array}{lllll}
\pa f\minidef 0, & 
\pa h \minidef gf \miniminus 1,   & 
\pa f_2 \minidef fh \miniminus lf, & 
\pa f_3 \minidef gf_2 \miniminus hh \miniplus g_2f,  & 
\pa f_4 \minidef ff_3 \miniminus lf_2 \miniplus f_2h \miniminus g_3f, ....
\\
\pa g \minidef 0, &
\pa l \minidef fg \miniminus 1,    &
\pa g_2 \minidef gl \miniminus hg, & 
\pa g_3 \minidef fg_2 \miniminus ll \miniplus f_2g, &
\pa g_4 \minidef gg_3 \miniminus hg_2 \miniplus g_2l \miniminus f_3g, ...
\end{array}
\]

A proof that the map $\alpha$ of~(\ref{honest}) is a homology
isomorphism is given in~\cite{markl:ip}.
We denote $\P_{(S,\underline \emptyset)}$ by $\Iso$
and resolution~(\ref{honest}) by $\alphaiso : \Riso \to \Iso$.
}
\end{example}

The above resolution was studied in great detail in~\cite{markl:ip} where
we called algebras over $\Riso$ {\em strong homotopy
equivalences\/}. It is not true that each (ordinary) chain homotopy
equivalence extends to a strong one. In~\cite{markl:ip} we, however, proved:

\begin{proposition}
\label{Dve_prasatka}
Let  $f,g,h,l$ be a chain homotopy equivalence as above. 
By changing either $h$ or
$l$, we can achieve that $(f,g,h,l)$ extends to a strong homotopy
equivalence, that is, to an algebra over the operad $\Riso$.
\end{proposition}

\noindent 
{\bf Side conditions in perturbation theory.}
Let us show how strong homotopy equivalences resonate with the classical
perturbation theory. The following material, as well as
Proposition~\ref{why_side_conditions}, will not be needed in the rest
of the paper.

Recall (see, for example, \cite{johansson-lambe:96}) that the
central object of perturbation theory, {\em strong deformation retract
(SDR)\/} data, are given by the diagram

\[
\label{SDR}
\SDR M{(A,\phi)}{\nabla}f
\]
\vglue4mm

\noindent 
of chain complexes and their maps, where $f\nabla  = \id_M$ and $\phi$
is a chain homotopy between $\nabla f$ and $\id_A$. One then shows
that one may always change the SDR-data in such
a way that the following {\em side conditions\/} hold:
\begin{equation}
\label{side_conditions}
\phi \phi=0,\ \phi \nabla = 0 \mbox { and } f \phi =0.
\end{equation}

The nature and r\^ole of these side conditions is explained in the following
proposition, proved in~\cite{markl:ip}. 

\begin{proposition}
\label{why_side_conditions}
Suppose we have a `partial' action $a: (\Gamma(f,g,h,l),\pa) \to
\End_{V,W}$ given by the SDR data as
\[
 a(f):= \nabla,\ a(g) := f,  a(h):= 0  \mbox { and } a(l) := \phi
\]
(we believe that the collision between $f$ denoting the generator and
$f$ denoting the map will not cause problems). Then $a$ can be
extended to a full action
\[
a: \Riso = (\Gamma(f,g,h,l,f_2,g_2,f_3,g_3,\ldots),\pa) \to \End_{V,W}
\]
by putting 
\begin{equation}
\label{Dan}
a(f_i) = a(g_i) = 0\mbox{ for } i \geq 2
\end{equation}
if and only if the side
conditions~(\ref{side_conditions}) are satisfied.
\end{proposition}

The difference between ordinary homotopy diagrams and strongly homotopy
ones is measured by the homology of the nerve of the category $S$. 
{}From this point of view, diagram~(\ref{hvezdicka}) related to
moves {\bf (M1)} and {\bf (M3)} was very complicated. 
On the other hand, homotopy versions of
diagrams related to Salvatore's move {\bf (S)} and move {\bf (M4)} are 
the obvious ones.

We close this section by a technical statement explaining how 
the `initial data' of the moves induce an
action of a cofibrant resolution. 
It is immediate to see that the operad $\P_\CD$ of
Definition~\ref{Krtecek} is the push-out of the diagram
\begin{equation}
\label{Pu}
\P_S \stackrel {n_1}{\vlla} \P_D \stackrel {n_2}{\vlra} \P_\uD
\end{equation}
where $\P_S := \P_{(S,\underline \emptyset)}$, $\P_D :=
\P_{(D,\underline \emptyset)}$, $n_1$ is the `corestriction' map and 
$n_2 := \phi_{(D,\underline \emptyset)}$. 
Let us assume that the operad $\P_D$ is free, $\P_D =
\Gamma(X_D)$, and that the maps $n_1$, $n_2$ are monomorphisms. These
conditions are satisfied for all moves considered in the paper.

Suppose also that $\alpha_S : \R_S \to \P_S$ (resp.~$\alpha_\uD :
\R_\uD \to \P_\uD$) is a cofibrant model for $\P_S$ (resp.~$\P_\uD$)
of the form $\R_S = (\Gamma(X^1_D,Y_S),d)$ (resp.~$R_\uD =
(\Gamma(X^2_D,Y_\uD),d)$), where $X^1_D$ (resp.~$X_D^2$) is  an identical copy
of $X_D$, $d(X^1_D) =0$ (resp.~$d(X^2_D) =0$) and $n_1 = \alpha_S \circ
\iota_1$ (resp.~$n_2 = \alpha_\uD \circ \iota_2$) , where $\iota_1$
(resp.~$\iota_2$) is the inclusion induced by the identification
$X^1_D \cong X_D$ (resp.~$X^2_D = X_D$). Such cofibrant models exist
by Theorem~\ref{factorization}. We will write, for $x_D \in X_D$,
$x^1_D := \iota_1(x_D)$ and $ x^2_D := \iota_2(x_D)$. 

\begin{lemma}
\label{sroubovacky}
Suppose that
\begin{equation}
\label{jsem_nachlazeny}
\alpha_S * \alpha_\uD : \calR_S * \calR_\uD \to \calP_S * \calP_D
\end{equation}
is a homology isomorphism. Then 
any two actions $A_S : \R_S \to \End_\bfU$ and $A_\uD : \R_\uD \to
\End_\bfU$ satisfying
\begin{equation}
\label{Bello}
A_S(x^1_D) = A_\uD (x^2_D) \mbox { for each $x_D \in X_D$}
\end{equation}
extend to an action
\[
A_{(S,\uD)} : \R_\CD \to \End_\bfU
\]
of a cofibrant resolution $\R_\CD$ of the operad $\P_\CD$.
\end{lemma}

\noindent
{\bf Remark.}
Intuitively, the above lemma means that a 
strongly homotopy $(S,\underline
\emptyset)$-diagram $A_S :\R_S \to \End_\bfU$ and a strongly homotopy
$\uD$-diagram $A_\uD : \R_\uD \to \End_\bfU$ that coincide on $D$ induce
a strongly homotopy $\CD$-diagram $A_{(S,\uD)} : \R_\CD \to \End_\bfU$
that extends both $A_S$ and $A_\uD$.

\noindent
{\bf Proof of the lemma.}
Let us denote by $\PU$ the push-out of~(\ref{Pu}), with the canonical
projection $\pi : \P_S * \P_\uD \to \PU$. Let us consider the diagram
\begin{center}
{
\unitlength=1.2pt
\begin{picture}(180.00,80.00)(0.00,0.00)

\put(160.00,50.00){\vector(3,2){0}}
\put(180.00,60.00){\makebox(0.00,0.00){$\End_\bfU$}}

\put(40.00,70.00){\makebox(0.00,0.00){%
     \scriptsize \hskip 10pt $\alpha_S * \alpha_\uD$}}
\put(120.00,70.00){\makebox(0.00,0.00){\scriptsize $\hskip 42pt A_S * A_\uD$}}
\put(150.00,20.00){\makebox(0.00,0.00){\scriptsize $A_\CD$}}
\put(40.00,10.00){\makebox(0.00,0.00){\scriptsize \hskip 15pt $\alpha_\CD$}}

\put(90.00,10.00){\vector(0,-1){0}}
\put(90.00,0.00){\makebox(0.00,0.00){$\R_\CD$}}
\multiput(70,00)(-2,0){20}{\makebox(0,0){$\cdot$}}
\multiput(100,5)(2,1.5){30}{\makebox(0,0){$\cdot$}}
\put(97.00,30.00){\makebox(0.00,0.00){\scriptsize $\iota$}}
\put(-10,30.00){\makebox(0.00,0.00){\scriptsize$\pi$}}

\put(90.00,60.00){\makebox(0.00,0.00){$\R_S * \R_\uD$}}
\multiput(90,50)(0,-2){20}{\makebox(0,0){$\cdot$}}

\put(30.00,0.00){\vector(-1,0){0}}
\put(0.00,0.00){\makebox(0.00,0.00){$\PU$}}

\put(0.00,60.00){\makebox(0.00,0.00){$\P_S*\P_\uD$}}
\put(0.00,50.00){\vector(0,-1){40.00}}
\put(110.00,60.00){\vector(1,0){50.00}}
\put(70.00,60.00){\vector(-1,0){50.00}}
\end{picture}}

\end{center}
in which $\alpha_S * \alpha_\uD$ is a homology isomorphism, 
by assumption~(\ref{jsem_nachlazeny}). We will construct a cofibrant resolution
$\alpha_\CD : \R_\CD
\to \PU$ by attaching free generators to
$\R_S*\R_\uD$. We will simultaneously construct also 
an action $A_\CD : \R_\CD \to
\End_\bfU$ extending $A_S * A_\uS$. These new data are indicated in the
above diagram by dotted lines.

The resolution $\R_\CD$ is constructed in the standard manner by killing
cycles in $\R_S * \R_\uD = (\Gamma(X^1_D,X^2_D,Y_S,Y_\uD),d)$ that
represent the homology of the kernel of $\pi(\alpha_S * \alpha_\uD)$
(compare the proof of Theorem~\ref{factorization}). 
The extension $A_\CD$ is
then given by an elementary obstruction theory.

We start by adding generators $\oX_D := \hskip 1mm \uparrow
\hskip -1mm X_D$ 
(the suspension of $X_D$) and defining the differential by
\[
d \overline x_D := x^1_D - x^2_D,\ x_D \in X_D.
\]
The extensions $\alpha_\CD$ and $A_\CD$ are given by
\begin{equation}
\label{nigra_sum}
\alpha_\CD (\overline x_D):= 0 \mbox { and }
A_\CD (\overline x_D):=0, \mbox { for $x_D \in X_D$.}
\end{equation}
Let us verify that these extensions commute with the differentials. For
$\alpha_\CD$ this means that $0 = \alpha_\CD(d \overline x_D)$, for
$x_D \in X_D$. This follows from
\[
\alpha_\CD(d \overline x_D) = \alpha_\CD(x^1_D - x^2_D) =
\alpha_\CD(x^1_D) - \alpha_\CD(x^2_D) = \pi(\alpha_S(x^1_S)) -
\pi(\alpha_\uD(x^2_D))
\]
since clearly $\pi(\alpha_S(x^1_S)) =
\pi(\alpha_\uD(x^2_D))$.
Similarly, 
\[
A_\CD(d \overline x_D) = A_\CD(x^1_D - x^2_D) =
A_\CD(x^1_D) - A_\CD(x^2_D) =0
\]
by~(\ref{Bello}), thus $A_\CD$ commutes with the differentials as well.

Since $\alpha_S * \alpha_\uD$ is a homology
isomorphism and ${\rm Ker}(\pi)$ is the ideal
generated in $\P_S * \P_\uD$ by $\{n_1(x_D)- n_2(x_D);\ x_D \in X_D\}$, 
each cycle in the kernel of
$\alpha_\CD : \Gamma(X^1_D,X^2_D,Y_S,Y_\uD,\oX_D) \to \PU$ 
is homologous to a cycle in the ideal $(\oX_D)$ generated
by $\oX_D$.  
Thus the `higher' generators $z$
killing
the cycles in the kernel of $\alpha_\CD$ 
can be added in such a way that $dz \in
(\oX_D)$. The extensions
\[
\alpha_\CD (z) :=  0\mbox { and } A_\CD(z) := 0
\]
then commute, by~(\ref{nigra_sum}), with the differentials.%
\qed

The following theorem verifies assumption~(\ref{jsem_nachlazeny}) in
two important cases.

\begin{theorem}
\label{takhle_prosvihnout_seky}
Let $\calP$ be an ordinary operad with $\calP(0)=0$, 
$\calP(1) \cong \bfk$ and trivial
differential. Let $\calS$ be either
$\calP$ or $\calP_\btb$ and let
$\alpha_\calS : \calR \to \calS$ be a cofibrant resolution.
Then the natural map 
\[
\alpha:= \alpha_\calS * \alphaiso : \calR * \Riso \to \calS * \Iso,
\]
where $\alphaiso: \Riso \to \Iso$ is the
resolution introduced in Example~\ref{tyden_pred_mistrovstvim}, is 
a homology isomorphism.
\end{theorem}

Observe that this theorem does not follow from
Theorem~\ref{easy_kunneth}, because neither $\Iso$ nor $\Riso$ are
augmented dg-operads. The proof is postponed to the Appendix.

\section{Proofs of the moves}
\label{jeste_jeden_medvidek}
 
\begin{proposition}[M1]
\label{Chrochtatko_a_Rypacek}
For each strongly homotopy $\P$-algebra $\bfV$, a chain complex
$W = (W,\partial)$ and a chain homotopy equivalence 
$f: (V,\pa) \to (W,\pa)$, there exist a strongly homotopy 
$\P$-structure $\bfW$ on 
$(W,\pa)$ and a strongly homotopy 
$\P$-structure ${\bf f}:\bfV \to \bfW$
on the map $f$.
\end{proposition}

\noindent
{\bf Proof.}
Let $\alpha_\bullet : \calR_\bullet \to \calP_\bullet$ be a cofibrant
resolution and $(S,D)$ the couple of diagrams introduced in Example~\ref{7}.
The initial data 
are given by an action $A_{\bullet \to \circ} : \R_\bullet*\Gamma(f) 
\to \End_{V,W}$
representing a strongly homotopy $\P$-algebra structure on $(V,\partial)$
and a chain map $f :(V,\pa) \to (W,\pa)$. 
Since $f$ is a chain homotopy equivalence
we can, by Proposition~\ref{Dve_prasatka}, extend $A_{\bullet \to
\circ}$, properly
choosing a homotopy inverse $g$ and homotopies $h$ and $l$, to a map
$A_\CD : \R_\bullet * \Riso \to \End_{V,W}$.

Very crucially, $\R_\bullet * \Riso$ is a cofibrant resolution of the
operad $\P_\CD$, because the operad $\P_\CD$ 
of~(\ref{Hopsalek}) is naturally isomorphic to
the free product $\P_\bullet * \Iso$ and the induced map
\begin{equation}
\label{Kremilek}
\alpha_\bullet * \alphaiso : 
\R_\bullet *\Riso \to \P_\bullet * \Iso \cong \P_\CD
\end{equation}
is, by Theorem~\ref{takhle_prosvihnout_seky}, 
a homology isomorphism. The above
situation is that of Lemma~\ref{sroubovacky} with $X_D = 0$.

\newcommand{\sssquare}[8]{
\setlength{\unitlength}{.7cm}
\begin{picture}(5,2)(0,2.5)
\thinlines

\put(0,4){\makebox(0,0){$#1$}}
\put(5,4){\makebox(0,0){$#2$}}
\put(0,1){\makebox(0,0){$#3$}}
\put(5,1){\makebox(0,0){$#4$}}

\put(-.5,2.5){\makebox(0,0)[r]{\scriptsize $#6$}}
\put(5.5,2.5){\makebox(0,0)[l]{\scriptsize $#7$}}
\put(2.5,1.5){\makebox(0,0)[b]{\scriptsize $#8$}}
\put(2.5,4.5){\makebox(0,0)[b]{\scriptsize $#5$}}

\put(1,1){\vector(1,0){3}}
\put(1,4){\vector(1,0){3}}
\put(0,3.5){\vector(0,-1){2}}
\put(5,3.5){\vector(0,-1){2}}
\end{picture}
}

By Theorem~\ref{factorization}, 
there exists a cofibration $\iota_\CD :\R_\bullet
*\Riso  \to \R_\uS$ and a homology isomorphism
$\alpha_\uS: \R_{\uS} \to
\P_{\uS}$ making the diagram
\begin{equation}
\label{Bara}
\sssquare{\R_\bullet *\Riso \hskip 2mm}%
       {\R_{\uS}}{\P_\bullet * \Iso}{\P_{\uS}}%
       {\iota_\CD}{\alpha_\bullet * \alpha_{\rm iso}}%
       {\alpha_{\uS}}{\phi_\CD}
\end{equation}

\vglue 1.2cm
\noindent 
commutative. Since we know, from Example~\ref{7}, 
that the couple $(D,S)$ has
the classical $\P$-extension property, Theorem~\ref{main} 
gives an action
$A_\uS: \R_{\uS} \to \End_{V,W}$.

We must prove the existence of strongly homotopy $\P$-structures 
on the complexes
$(V,\pa)$ and $(W,\pa)$ and a sh{} structure on the map $f$, i.e.~of
an action $A_\btb:\R_{\btb} \to \End_{V,W}$, 
where $\R_\btb$ is a cofibrant
resolution of the operad $\P_\btb$ of~(\ref{dve_prasatka}). 
We also require that $A_{\btb}$
extends the given data, i.e.~that there is a 
cofibration $\iota_{\bullet \to \circ} : \P_\bullet *
\Gamma(f) \to \R_\btb$ such that the diagram

\begin{equation}
\label{B}
\sssquare{\R_\bullet* \Gamma(f) \hskip 4mm}{\R_\btb}%
          {\P_\bullet* \Gamma(f) \hskip 3mm}{\P_\btb}%
          {\iota_{\bullet \to \circ}}%
          {\alpha_{\bullet \to \circ} := \alpha_\bullet*\id}%
          {\alpha_\btb}%
          {\phi_{\bullet \to \circ}}
\end{equation}

\vskip 1.2cm
\noindent 
commutes and that $A_\btb (\iota_{\bullet \to \circ}) 
= A_{{\bullet \to \circ}}$.
Consider the diagram

\begin{center}
{
\unitlength=1.3pt
\begin{picture}(220.00,115.00)(0.00,0.00)
\put(150.00,10.00){\makebox(0.00,0.00){\scriptsize $r$}}
\put(50.00,10.00){\makebox(0.00,0.00){\scriptsize $\alpha_\btb$}}
\put(210.00,0.00){\makebox(0.00,0.00){$\P_\uS$}}
\put(110.00,0.00){\makebox(0.00,0.00){$\P_\btb$}}
\put(10.00,0.00){\makebox(0.00,0.00){$\R_\btb$}}
\put(130.00,20.00){\makebox(0.00,0.00){\fbox{$4$}}}
\put(167.00,30.00){\makebox(0.00,0.00){\scriptsize $ \phi_\CD$}}
\put(93.00,30.00){\makebox(0.00,0.00){\scriptsize $ \phi_{\bullet \to \circ}$}}
\put(35.00,30.00){\makebox(0.00,0.00){\fbox{$1$}}}
\put(180.00,70.00){\makebox(0.00,0.00){\fbox{$2$}}}
\put(220.00,50.00){\makebox(0.00,0.00){\scriptsize $ \alpha_\uS$}}
\put(160.00,50.00){\makebox(0.00,0.00){$\P_\bullet * \Iso$}}
\put(100.00,60.00){\makebox(0.00,0.00){\scriptsize $ \id * j$}}
\put(60.00,50.00){\makebox(0.00,0.00){$\P_\bullet * \Gamma(f)$}}
\put(0.00,50.00){\makebox(0.00,0.00){\scriptsize $ \iota_{\bullet \to \circ}$}}
\put(150.00,80.00){\makebox(0.00,0.00){%
                                 \scriptsize $ \alpha_\bullet * \alphaiso$}}
\put(90.00,80.00){\makebox(0.00,0.00){\fbox {$3$}}}
\put(50.00,80.00){\makebox(0.00,0.00){\scriptsize $\alpha_\bullet * \id$}}
\put(210.00,100.00){\makebox(0.00,0.00){$\R_\uS$}}
\put(160.00,110.00){\makebox(0.00,0.00){\scriptsize $ \iota_\CD$}}
\put(110.00,100.00){\makebox(0.00,0.00){$\R_\bullet * \Riso$}}
\put(50.00,110.00){\makebox(0.00,0.00){\scriptsize $ \id * \tilde \j{}$}}
\put(10.00,100.00){\makebox(0.00,0.00){$\R_\bullet *\Gamma(f)$}}

\put(210.00,90.00){\vector(0,-1){80.00}}
\put(10.00,90.00){\vector(0,-1){80.00}}
\put(20.00,90.00){\vector(1,-1){30.00}}
\put(170.00,40.00){\vector(1,-1){30.00}}
\put(120.00,90.00){\vector(1,-1){30.00}}
\put(70.00,40.00){\vector(1,-1){30.00}}

\put(23.00,0.00){\vector(1,0){74.00}}
\put(123.00,0.00){\vector(1,0){77.00}}

\put(85.00,50.00){\vector(1,0){55.00}}

\put(130.00,100.00){\vector(1,0){70.00}}
\put(30.00,100.00){\vector(1,0){60.00}}

\end{picture}}
\end{center}
where $j: \Gamma(f) \to \Iso$ maps $f$ to the corresponding element of
$\Iso$, $\tilde {\j{}} : \Gamma(f) \to \Riso$ is a lift of $j$, i.e.~a map
that satisfies $\alpha_{{\rm iso}} (\tilde \j) = j$ (such a map
clearly exists) and $r : \P_{\bullet \to \bullet} \to \P_\uS$ is the
corestriction map.

The subdiagram \fbox1 (resp.~\fbox2\hskip.5mm)
commutes because it is~(\ref{B})
(resp.~(\ref{Bara})). The commutativity of
the remaining subdiagrams \fbox3 and \fbox4 is obvious. This implies
the commutativity of the `outer diagram'
 
\begin{center}
\square{\R_\bullet* \Gamma(f)\hskip 6mm}{\R_{\uS}}{\R_{\bullet \to \bullet}}%
       {\P_{\uS}}{\iota_\CD (\id * \tilde \j{})}
       {\iota_{\bullet \to \circ}}{\alpha_{\uS}}
       {r(\alpha_{\bullet \to \bullet})}
\end{center}
By~Lemma~\ref{ll} there exists a map $\beta:
\R_{\bullet \to \bullet} \to \R_\uS$ such that $\beta(
\iota_{\bullet \to \circ}) = \iota_\CD(\id * \tilde \j)$. The action
$A_{\bullet \to \bullet} := A_\uS (\beta) : \R_\btb \to \End_{V,W}$ 
then has the required property.
\qed

\begin{proposition}[M2]
Suppose $\bfV$ and $\bfW$ are two strongly
homotopy $\P$-algebras and ${\bf f}:\bfV \to
\bfW$ a strongly homotopy  
$\P$-algebra morphism. Suppose that $g:(V,\pa) \to (W,\pa)$
is a chain map that is chain homotopic to
$\bff_{\#}$. Then there exists a strongly
homotopy $\P$-structure on $g$.
\end{proposition}

\noindent
{\bf Proof of (M2)}.
Let $f,g:  (U,\pa) \to (V,\pa)$ 
be chain maps and $h$ a chain homotopy between $f$ and
$g$. 
Let $\alpha_\btb :\Rbtb \to \calP_\btb$ be a cofibrant resolution
of the operad $\calP_\btb$ introduced in~(\ref{dve_prasatka}). 
Let $\#: \Gamma(u) \to
\Rbtb$ be as in~(\ref{dnes_jsem_blue}). Finally, 
let $\bff : \calR_\btb \to \End_{U,V}$ be a sh{}
$\calP$-structure on the dg map $f$, that is, $\bff(\#u) = f$ 
(cf.~Definition~\ref{blue_mood}). 

Consider the free acyclic bi-colored operad $\Gammax$ on generators
$\xi,\pa\xi : \bullet \to \bullet$ of degrees $1$ and $0$,
respectively. By Theorem~\ref{easy_kunneth}, the homomorphism
$\beta :\Rbtb * \Gammax \to \calP_\btb$ defined by $\beta|_\Rbtb :=
\alpha_\btb$ and $\beta(\xi) = \beta(\pa\xi) := 0$,
is also a cofibrant resolution of $\Pbtb$. Now introduce
\[
B : \Rbtb * \Gammax \to \End_{U,V}
\]
by $B|_{\Rbtb} := \bff$, $B(\xi) := h$ and $B(\pa\xi) := g-f$. 
Finally, define $\jmath : \Gamma(u) \to \Rbtb * \Gammax$ as 
$\jmath(u) := \#(u) +
\pa\xi$.

Lemma~\ref{ll} gives a map $\gamma :
\Rbtb \to \Rbtb * \Gammax$ in the diagram:
\begin{center}
\squarewithdotedmodified{\Gamma(u)}{\Rbtb*\Gammax}{\Rbtb}{\calP_\btb}%
       {\jmath}{\#}{\beta}{\alpha_\btb}{\gamma} 
\end{center}
such that $\gamma \circ \# = \jmath$. Since then
\[
(B\circ \gamma)(\# (u)) = B(\jmath(n)) = B(\#(u) + \pa\xi) = f +(g-f) =g, 
\]
${\bf g}:= B\circ \gamma :\Rbtb \to \End_{U,V}$ 
is the desired strongly homotopy $\calP$-structure on $g$.%
\qed

\begin{proposition}[M3]
\label{Das ist Blut}
Suppose ${\bf f}:\bfV\to \bfW$ is an sh{}
$\P$-algebra map such that $\bff_{\#}$ is a chain homotopy
equivalence. 
Suppose that  $g:(W,\pa) \to (V,\pa)$ is a chain homotopy
inverse of $\bff_{\#}$. Then there exists a strongly
homotopy $\P$-structure on the map $g$.
\end{proposition}

\noindent
{\bf Proof.}
The `initial data' consist of an action $A_\btb : \R_\btb \to
\End_{V,W}$ representing ${\bf f}
: {\bf V} \to {\bf W}$, and an action $A_{\rm iso} : \Riso \to \End_{V,W}$,
induced by $f := {\bf f}_{\#}$ and its homotopy inverse $g$.

The assumptions of Lemma~\ref{sroubovacky} are clearly
satisfied, with $S$ the diagram
\begin{center}
{
\unitlength=.45pt
\begin{picture}(180.00,40.00)(80.00,30.00)
\put(20.00,10.00){\vector(-3,1){1.00}}
\put(80.00,52.00){$f$}
\put(80.00,-22.00){$g$}
\put(160.00,30.00){\vector(3,-1){1.00}}
\bezier{100}(20.00,10.00)(90.00,-10.00)(160.00,10.00)
\bezier{100}(20.00,30.00)(90.00,50.00)(160.00,30.00)
\put(180.00,20.00){\makebox(0.00,0.00){$\circ$}}
\put(0.00,20.00){\makebox(0.00,0.00){$\circ$}}
\put(187.00,20.00){\makebox(0.00,0.00)[l]{, $fg=1$ and $gf =1$,}}
\end{picture}}
\end{center}
\vskip .8cm

\noindent 
$D:= \circ \hskip 2mm \lra \hskip 2mm \circ$, 
and $X_D := {\rm Span}(f)$; assumption~(\ref{jsem_nachlazeny}) is
guaranteed by Theorem~\ref{takhle_prosvihnout_seky}. Therefore the
initial data induce an action $A_\CD : \R_\CD \to \End_{V,W}$ of a
cofibrant resolution $\R_\CD$ of the operad $\P_\CD$. 

As we showed in Example~\ref{9}, the couple $(S,D)$ has the classical
extension property and Theorem~\ref{main} guarantees the existence of
an action $A_\uS : \R_\uS \to \P_\uS$ of a cofibrant resolution of the
operad $\P_\uS$.
The rest follows, as in the proof of Proposition~\ref{Chrochtatko_a_Rypacek}, 
from a general nonsense argument.
\qed

In exactly the same manner, one can prove {\bf (M4)}
and also Salvatore's move {\bf (S)}. The
situation is even easier than above, because the strongly homotopy
diagrams related to these moves are the ordinary homotopy disgrams and also
Lemma~\ref{sroubovacky} applies in a straightforward manner.

\section{Appendix (remaining proofs)}
\label{jeden_medvidek}

{\bf Proof of Theorem~\ref{factorization}}. 
As in the proof of~\cite[Theorem~2.1.1]{hinich:CA97}, 
we construct inductively a sequence of cofibrations 
$C_{-1} := \calS \hookrightarrow C_0 \hookrightarrow 
C_1 \hookrightarrow C_2 \hookrightarrow \cdots\ $ and extensions
$\alpha_n : C_n \to \calQ$ of $f: \calS \to \calQ$
such that $C_n = C_{n-1}*\Gamma(E_n)$ 
with $\pa (E_n) \subset C_{n-1}$,  $n \geq 0$, with the following properties:
\begin{itemize}
\item[(i)]
$\alpha_n : C_n \to \calQ$ is an epimorphism,
\item[(ii)]
$Z(\alpha_i) : Z(C_n) \to Z(\calQ)$, where $Z(-)$ denotes cycles, 
is also an epimorphism and
\item[(iii)]
whenever $z \in Z(C_n)$ is mapped by $\alpha_i$ to a boundary in $\calQ$,
then $z$ is a boundary in~$C_{n+1}$.
\end{itemize}
Observe that if the conditions (i) and (ii) are satisfied with $n=0$, then they are
satisfied with an arbitrary $n \geq 0$.

The first step of the construction can be made as follows.
Let $E_0 := \calQ \oplus\! \desusp \! \calQ$ with $\calQ$ interpreted as a
colored collection and $\desusp \! \calQ$ as its desuspension. The
differentials are defined by $\pa q := \desusp
q$, $\pa(\desusp q) := 0$, and the map $\alpha_0$ by 
$\alpha_0(q): = q$ and $\alpha_0(\desusp q) :=
\pa q$, $q \in \calQ$. This choice clearly fulfills (i) and (ii) with $n=0$.
Then we proceed by induction.

For $n \geq 1$, let $E_n$ be the collection spanned by all pairs 
$(z,u)$ with $z$ a cycle in $C_n$ and $u \in \calQ$ such that
$\alpha_{n-1}(z) = \pa u$. We extend the differential by $\pa (z,u):
=z$ and put $\alpha_i(z,n):= u$.

The operad $\calR := \dirlim C_n$ together with $\alpha := \dirlim
\alpha_n :\calR \to \calQ$ are objects whose existence we require.%
\qed

\noindent
{\bf Proof of Lemma~\ref{ll}.}
Let $\iota: (\calS,\pa) \hookrightarrow 
(\calS* \Gamma(E),\pa) = \calQ$ be as in Definition~\ref{cof}. 
Let $h_{-1} := f : \calS \to \calX$ and suppose we
have already defined, for some $n \geq 0$, 
\[
h_{n-1} : \calS * \Gamma(E)_{< n} \to \calX
\]
such that $h_{n-1}|_\calS = f$ and $p \circ h_{n-1} 
= g|_ {\calS * \Gamma(E)_{< n}}$. 
Let us extend $h_{n-1}$ to some $h_{n}: \calS * \Gamma(E)_{\leq n} \to
\calX$ 
with similar properties by showing that, for each $e
\in E_{n}$, there exists a solution $\omega$ of the equations
\[
\pa \omega = h_{n-1}(\pa e)
\ \mbox { and }\
p (\omega) = g(e);
\]
then $h_n(e) := \omega$ has the desired properties.
We find $\omega$ in the form  $\omega = a +
b$, where 
\begin{equation}
\label{polevka_mi_pomuze}
p(a) = g(e)
\end{equation}
and $b \in \Ker(b)$ satisfies
\[
\pa b = h_{n-1}(\pa e) - \pa a.
\]
A solution $a$ of~(\ref{polevka_mi_pomuze}) exists 
because $p$ is an epimorphism. Then $ h_{n-1}(\pa e) - \pa a \in \Ker(p)$ and
the existence of $b$ follows from the acyclicity of $\Ker(p)$ which in
turn follows from the long exact sequence for the homology isomorphism
$p: \calX \epi \calY$. The construction of $h$ is completed by induction.%
\qed

\noindent 
{\bf Proof of Theorem~\ref{easy_kunneth}.}
To understand better the idea of the proof and the r\^ole of
the augmentation, consider a similar
statement for dg-associative algebras. An associative
augmented algebra $A=(A,\pa)$ decomposes as $A \cong \bfk \oplus \overline A$,
where the augmentation ideal $\overline A$ is $\pa$-closed. If $B =
\bfk \op \overline B$ is
another augmented associative algebra, then clearly
\begin{equation}
\label{tag}
A*B \cong \bfk \ \op \ \oA \ \op \  \oB \ \op \  
\oA\ot \oB \ \op \  \oB \ot \oA \ \op \ 
\oA\ot \oB \ot \oA \ \op \  \oB \ot \oA \ot \oB \ \op \  \cdots . 
\end{equation}
The ordinary K\"unneth theorem~\cite{hilton-stammbach} then implies that
\begin{eqnarray*}
\lefteqn{
H(A*B) \cong
H(\bfk \ \op \ \oA \ \op \  \oB \ \op \  
\oA\ot \oB \ \op \  \oB \ot \oA \ \op \ 
\oA\ot \oB \ot \oA \ \op \  \oB \ot \oA \ot \oB \ \op \  \cdots)}
\\
&&
\hskip -3mm \cong 
 H(\bfk) \ \op \ H(\oA) \ \op \  H(\oB) \ \op \  
H(\oA)\ot H(\oB) \ \op \  H(\oB) \ot H(\oA) \ \op \ 
H(\oA)\ot H(\oB) \ot H(\oA) \ \op  \cdots  
\\
&&
\hskip -3mm\cong 
H(\bfk) \ \op \ \overline{H(A)} \ \op \  \overline{H(B)} \ \op \  
\overline{H(A)}\ot \overline{H(B)} \ \op \  \overline{H(B)} \ot 
\overline{H(A)} \ \op \ 
\overline{H(A)}\ot \overline{H(B)} \ot \overline{H(A)} \ \op  \cdots
\\
&&
\cong H(A)*H(B).
\end{eqnarray*}

Let us move to the operadic case.
If $\calP_1$ and $\calP_2$ are augmented operads with augmentation ideals
$\overline{\calP}_s$, $s=1,2$, then the free product $\calP_1 *
\calP_2$  has a decomposition similar to~(\ref{tag}), given by a
summation over isomorphism classes of trees
(see~\cite{markl-shnider-stasheff:book} for the terminology). 

Let us describe this decomposition for non-colored operads first.
We say that a rooted tree $T$ is vertex-labeled, if it is equipped with
a map $m : {\it Vert}(T) \to \{1,2\}$ from its set 
of internal vertices ${\it Vert}(T)$. Then the augmentation ideal of
the augmented operad $\calP_1 * \calP_2$ is given as
\begin{equation}
\label{zpatky_z_Izraele}
(\overline{\calP_1 * \calP_2})(n) =  \bigoplus_{T} \bigotimes_{v \in {\it
Vert(T)}} \overline{\calP}_{m(v)} (\# {\it Vert}(v)), 
\end{equation}
where the summation is taken over representatives of isomorphism
classes of vertex-labeled $n$-trees. 
The above representation is canonical up to
choices of representatives of isomorphism classes of trees, 
compare~\cite[Proposition~1.82]{markl-shnider-stasheff:book}.

For colored operads, there exist a decomposition similar
to~(\ref{zpatky_z_Izraele}), but based on vertex-labeled trees with 
{\bf C}-colored leaves.  Theorem~\ref{easy_kunneth} then follows 
from the ordinary K\"unneth theorem applied to~(\ref{zpatky_z_Izraele}).%
\qed

\noindent 
{\bf Proof of Theorem~\ref{takhle_prosvihnout_seky}}.
We may suppose, by Proposition~\ref{zitra_brzo_vstavam}, that 
$\calR = (\Gamma(X),\pa_\calR)$ with $X = X^0 \op X^1 \op
\cdots$ is bigraded.
Recall that the operad 
$\Riso = (\Gamma(f,g,h,l,f_2,g_2,\ldots), \paiso)$ was introduced
in Example~\ref{tyden_pred_mistrovstvim}.
Then 
\[
\calR * \Riso = (\Gamma(X,f,g,h,l,f_2,g_2,\ldots),\pa)
\] 
with the differential $\pa$ induced by $\paiso$ and $\parr$ in the standard
way. There is a natural `upper' grading of $\calR * \Riso$  
induced by TJ-grading of $\calR$
and the usual homological grading of $\Riso$. Differential $\pa$ has
degree $-1$ with respect to this grading.
A direct inspection shows that  $\alpha_\calS * \alphaiso$ induces an
isomorphism 
\[
H^0(\calR * \Riso) \cong \calS * \Iso
\]
therefore it suffices to prove that 
$\calR * \Riso$ is acyclic in
positive upper degrees, that is
\[
H^{> 0}(\Gamma(X,f,g,h,l,f_2,g_2,\ldots),\pa) = 0.
\] 
Observe that this statement does not follow
from Theorem~\ref{easy_kunneth}, because $\Riso$ is not a dg-augmented
operad.
We derive the  acyclicity from the following sublemma.

\begin{sublemma}
The inclusion 
$\Gamma(X,f,g) \hookrightarrow \Gamma(X,f,g,h,l,f_2,g_2, \ldots)$ 
is a homology epimorphism.
\end{sublemma}

Assuming the sublemma, $H^{> 0}(\Gamma(X,f,g),\pa) \to 
H^{> 0}(\Gamma(X,f,g,h,l,f_2,g_2\ldots),\pa)$ is an epimorphism. Since
$\pa|_{\Gamma(X,f,g)} = \pa_0$ is induced by $\pa_{\calR}$, 
the operad $(\Gamma(X,f,g),\pa)$ is dg-augmented and
Theorem~\ref{easy_kunneth} implies that 
\[
H^{> 0}(\Gamma(X,f,g),\pa) 
\cong H^{> 0}(\Gamma(X),\pa_{\calR}) * \Gamma(f,g) = 0,
\]
because $H^{>
0}(\Gamma(X),\pa_\calR)=0$ by (iv) of
Definition~\ref{zitra_jedu_do_Oberwolfachu}.

It remains to prove the sublemma. We need a third grading, 
the h-grading (h from homogeneity) of
$\Gamma(X,f,g,h,l,f_2,g_2, \ldots)$ given by the number of
generators $f,g,h,l,f_2,g_2,\ldots$ (so $X$ has h-grading $0$). 
With respect to this h-grading, differential $\pa$
decomposes to the sum
\[
\pa = \pa_{-1} + \pa_0 + \pa_{+1}
\]
where $\pa_{-1}$ and $\pa_{+1}$ are induced by $\paiso$ and
$\pa_0$ is induced by $\parr$. 
Recall that in~\cite{markl:ip} we proved that
$(\Gamma(f,g,h,l,f_2,g_2, \ldots), \pa_{+1})$ is acyclic in positive
homological degrees. Since $(\Gamma(f,g,h,l,f_2,g_2, \ldots), \pa_{+1})$ is
dg-augmented and
\[
(\Gamma(X,f,g,h,l,f_2,g_2, \ldots), \pa_{+1}) \cong (\Gamma(X),\pa=0) * 
(\Gamma(f,g,h,l,f_2,g_2, \ldots), \pa_{+1}),
\] 
we may use Theorem~\ref{easy_kunneth} to derive that 
the ideal ${\cal I}$ generated in 
$\Gamma(X,f,g,h,l,f_2,g_2, \ldots)$ by $h,l,f_2,g_2, \ldots$ is 
$\pa_{+1}$-acyclic. Since 
\[
\Gamma(X,f,g,h,l,f_2,g_2,\ldots) \cong \Gamma(X,f,g) \oplus \calI
\]
the sublemma will be proved if we show that each
$\pa$-cycle $a$ of $\Gamma(X,f,g,h,l,f_2,g_2,\ldots)$ is $\pa$-homologous to a
$\pa$-cycle of $\Gamma(X,f,g)$. 

Suppose that $a \not\in \Gamma(X,f,g)$, that is $a = x + y$, where 
$0 \not=x \in {\cal I}$ and $y \in \Gamma(X,f,g)$.
We may decompose $x$ as $x_1+\cdots+
x_N$ with respect to the h-degree. Then $\pa_{+1}(x_N) = 0$ because
this is precisely the component of $\pa(a)$ in $\calI$ of h-grading $N+1$.
Because $\calI$ is $\pa_{+1}$-acyclic, $x_N =
\pa_{+1}(b_{N-1})$ for some $b_{N-1} \in \calI$. 
The highest h-component of $x' := x - \pa(b_{N-1})$ has h-degree $N-1$.

We may apply the above considerations to $a' := a - \pa(b_{N-1})$ and
repeat this process $N-1$ times. We end up with a decomposition $a =
x+y$, where $x= x_1$ is of h-degree $1$. Since $\calI$ does not
contain elements of h-degree $0$, $\pa_+(x_1) = 0$ together with the
$\pa_{+1}$-acyclicity of $\calI$ implies
that $x_1 =0$, therefore $a \in \Gamma(X,f,g)$.
The proof of the sublemma is finished.%
\qed 

{\small
\def\cprime{$'$}

}

\vfill
\hfill
{\tt \jobname.tex}

\end{document}